\documentclass[11pt]{article}

\usepackage[T1]{fontenc}
\usepackage[utf8]{inputenc}
\usepackage{lmodern}
\usepackage{amsmath,amssymb,amsthm}
\usepackage[margin=2.6cm]{geometry}
\usepackage{booktabs}
\usepackage{array}
\usepackage{xcolor}
\usepackage[colorlinks=true,linkcolor=blue!50!black,citecolor=blue!50!black,urlcolor=blue!50!black]{hyperref}
\usepackage[nameinlink]{cleveref}

\makeatletter
\renewcommand{\@makefnmark}{\hbox{\@textsuperscript{\normalfont(\@thefnmark)}}}
\makeatother
\newtheorem{theorem}{Theorem}
\newtheorem{proposition}{Proposition}
\newtheorem{lemma}{Lemma}
\newtheorem{corollary}{Corollary}
\theoremstyle{definition}
\newtheorem{criterion}{Criterion}
\newtheorem{observation}{Observation}
\theoremstyle{remark}
\newtheorem{remark}{Remark}
\newcommand{\NP}{\boldsymbol{\mathcal{S}}}
\newcommand{\Mx}{\max_{y\in\mathbb{R}^3}\Big(\frac12|v|^2+\lambda\big(|\omega|^2+\varepsilon^2\big)^{1/2}\Big)(y,t)}
\newcommand{\Mq}{\max_{y\in\mathbb{R}^3}\big(|v|^4+\nu^2|\omega|^2\big)(y,t)}

\numberwithin{equation}{section}

\title{Local Reynolds number criteria for the regularity\\
of the three-dimensional Navier-Stokes equations,\\
and their consistency with forced finite-time blow-up}
\author{Ali Bouzekri\\[2pt]
\small Laboratoire de Math\'ematiques Pures et Appliqu\'ees (LMPA), Tizi-Ouzou, Algeria\\
\small\texttt{bouzekri.ali.Lmpa@gmail.com}}
\date{November 2026}

\begin{document}
\maketitle

\begin{abstract}
{We give four regularity criteria for the forced Navier-Stokes equations in $\mathbb{R}^3$. The first is an $\dot H^1$ estimate: a strong solution cannot blow up while the local Reynolds number $|v|^2/(\nu|\omega|)$ stays below $4/(3\sqrt3)$ wherever the vorticity exceeds a threshold $M(t)$ with $\int M^{6/5}\,dt<\infty$. It rests on an exact identity that gives the cubic vorticity term a favorable sign. Two further $\dot H^1$ estimates use the Newton potential of $|\omega|^2$, either as a multiplier or as a gradient added to the Lamb vector, so that the integral of $|\omega|^2|v|^2$ cancels. The second of these gives a critical condition on the vorticity alone, implied by $\int\|\omega\|_{L^3}^2\,dt<\infty$. The fourth is a maximum principle for $\frac12|v|^2+\lambda|\omega|$: vortex stretching is absorbed by viscous dissipation, and the pressure is controlled by a weighted local estimate. It requires a bound on $|v|^2/(\nu|\nabla v|)$ only at maximum points where the velocity part dominates. We show that the forced finite-time blow-up recently constructed by OpenAI violates each hypothesis, and that the violations locate its mechanism: inflow down the centripetal pressure gradient, at the same order as viscous dissipation.}
\end{abstract}

\section{Introduction}

We consider the incompressible Navier-Stokes equations in $\mathbb{R}^3$ with an external force,
\begin{equation}\label{eq:NS}
\partial_t v+(v\cdot\nabla)v-\nu\Delta v+\nabla p=f,\qquad \nabla\cdot v=0,\qquad v(\cdot,0)=v_0,
\end{equation}
with $\nu>0$. The vorticity is $\omega=\nabla\times v$, $\langle\cdot,\cdot\rangle$ denotes the pairing of $L^2(\mathbb{R}^3)$, matrix norms are Frobenius norms \footnote{That is for a matrix \(M\) with entries \((a_{ij})_{ij}\) \(|M|=\sqrt{\sum a_{ij}^2}\)}, and $S=\left((\nabla v)+(\nabla v)^T\right)/2$ denotes the symmetric part of $\nabla v$. Since the antisymmetric part of $\nabla v$ has squared norm $|\left((\nabla v)-(\nabla v)^T\right)/2|^2=\frac12|\omega|^2$,
\begin{equation}\label{eq:split}
|\nabla v|^2=|S|^2+\frac12|\omega|^2 .
\end{equation}

The quantities compared in this paper are $|v|^2$ and $\nu|\omega|$, or $|v|^2$ and $\nu|\nabla v|$. Their ratios are local Reynolds numbers: the local speed $|v|$ times the length $|v|/|\omega|$ (or $|v|/|\nabla v|$) over which it varies, divided by $\nu$. They are invariant under the scaling that preserves~\eqref{eq:NS}, and \Cref{sec:scaling} shows that no change of variables, including one that changes the viscosity, alters a bound on them. A hypothesis bounding a local Reynolds number is therefore critical, and its constant is a pure number.

\paragraph{Results.} \Cref{thm:main} proves regularity when $|v|^2\le C\nu|\omega|$ with $C<4/(3\sqrt3)$ on the set where $|\omega|>M(t)$, with $\int M^{6/5}<\infty$. \Cref{thm:newton} replaces the multiplier $|\omega|v$ of that proof by $\NP(|\omega|^2)v$, with $\NP$ the inverse of $-\Delta$. The integral of $|\omega|^2|v|^2$ then cancels identically, the only remaining term is $\int\NP(|\omega|^2)^2|v|^2$, and regularity follows when the time integral of $\int\NP(|\omega|^2)^2|v|^2\,dx\big/\big(\nu^3\|\nabla v\|_{L^2}^2\big)$ is finite, a critical condition without smallness, implied for instance by $\int\|\omega\|_{L^{12/5}}^2\|\omega\|_{L^2}\|v\|_{L^3}\,dt<\infty$. \Cref{thm:coulomb} adds instead the gradient $\nabla\NP(|\omega|^2)$ to $\omega\times v$. This gradient is orthogonal to $\Delta v$, the same cancellation occurs, and regularity follows when the time integral of $\int\NP(|\omega|^2)|\omega|^2\,dx\big/\big(\nu\|\nabla v\|_{L^2}^2\big)$ is finite, a critical condition on the vorticity alone, implied by $\int\|\omega\|_{L^3}^2\,dt<\infty$. \Cref{thm:max} is a maximum principle for $\frac12|v|^2+\lambda|\omega|$, with $\lambda$ a constant proportional to $\nu$. It proves regularity when $|v|^2\le C\nu|\nabla v|$ with $C$ below an explicit universal constant, at the near-maximum points where the velocity part dominates. Its proof uses a weighted estimate of $\nabla v$ by $\omega$ (\Cref{lem:weighted}) that controls the pressure gradient at a point by $\|v\|_{L^\infty}\|\omega\|_{L^\infty}$ (\Cref{cor:pressure}). \Cref{prop:signed} records the elementary maximum principle under a hypothesis on the signed pressure work. \Cref{sec:forced} and \Cref{rem:newton-forced}, \ref{rem:coulomb-forced} and~\ref{rem:lamb-forced} check every hypothesis against the forced blow-up of~\cite{OpenAI26}.
\paragraph{Relation to known criteria.} Regularity criteria involving velocity and vorticity jointly are numerous. Closest to \Cref{thm:main} are local criteria for suitable weak solutions placing the component of the velocity orthogonal to the vorticity direction, $v\times\omega/|\omega|$, or the component of the vorticity orthogonal to the velocity direction, $\omega\times v/|v|$, in scale-invariant Lebesgue classes on a parabolic cylinder~\cite{Chae07}. Chae and Lee later improved them to a criterion on the positive part of the triple product of $v\times\omega/|\omega|$ with the direction of $\nabla\times\omega=-\Delta v$~\cite{CL16}, which is the triple product appearing in the $\dot H^1$ estimate below. Miao and Wang obtained criteria of the same type with small local scaled norms from the rotation form of the equations, together with criteria on the Bernoulli pressure $\frac12|v|^2+p$~\cite{MW17}. A further criterion constrains the vorticity direction only where the vorticity is large~\cite{LRT25}. These are integrability or direction conditions, and the local ones are proved for suitable weak solutions through the local energy inequality and $\varepsilon$-regularity on parabolic cylinders. Our proofs use no local regularity theory. \Cref{thm:main}, \ref{thm:newton} and~\ref{thm:coulomb} are global $\dot H^1$ estimates closed by Gronwall's lemma, and \Cref{thm:max} is a maximum principle with an explicit pointwise bound on the pressure gradient; the only exception is \Cref{cor:lamb-single}, which relies on the small-data theory in $BMO^{-1}$. The price is that our criteria are global: they constrain the whole flow, not a neighbourhood of a point. The hypotheses of \Cref{thm:main} and~\ref{thm:max} are pointwise comparisons of magnitudes with explicit constants, and we have found no reduction of them to those conditions. The most direct attempt fails: the hypothesis of \Cref{thm:main} gives $|v\times\omega/|\omega||\le(C\nu)^{1/2}|\omega|^{1/2}$, and placing the right-hand side in the scale-invariant class required in~\cite{Chae07} amounts to the critical vorticity condition $\frac2q+\frac3p=2$ of~\cite{BdV95}, which is not available a priori. Among pressure criteria, regularity holds under a lower bound on the pressure or on the Bernoulli quantity $|v|^2+2p$~\cite{SS02}. The observation that a velocity peak can grow only where the pressure work $-v\cdot\nabla p$ is positive underlies the criteria of~\cite{TY16}, which are localised to the high-velocity region. \Cref{prop:signed} is the pointwise form of that observation; we include it because it is the criterion closest to the mechanism of the forced blow-up. \Cref{thm:max} differs from these: its hypothesis involves only the velocity and its gradient, and the pressure is controlled through \Cref{lem:weighted}.

\paragraph{Relation to the forced blow-up.} A finite-time blow-up for~\eqref{eq:NS} with zero initial data and a smooth, compactly supported force was announced in~\cite{OpenAI26}; see also~\cite{Dur26,LR26} for independent readings of the construction and~\cite{BCH26} for a discussion of its formalisation. It addresses the forced alternatives of the Millennium Problem formulation~\cite{Fef06}; the unforced alternatives remain open. The data of that solution satisfy every hypothesis on $v_0$ and $f$ made here, so every regularity criterion must fail on it through its own hypothesis. \Cref{sec:forced} and \Cref{rem:newton-forced}, \ref{rem:coulomb-forced} and~\ref{rem:lamb-forced} check this for each criterion and measure by how much each fails.

\section{Scaling and the local Reynolds number}\label{sec:scaling}

 For $(v,p,f)$ solutions to~\eqref{eq:NS} with viscosity $\nu$ on $\mathbb{R}^3\times[0,T)$, and let $a,b>0$. We have the following invariance
\begin{equation}\label{eq:scaling}
v'(x,t)=a\,v(bx,abt),\qquad p'(x,t)=a^2p(bx,abt),\qquad f'(x,t)=a^2b\,f(bx,abt)
\end{equation}
solve~\eqref{eq:NS} with viscosity $\nu'=\nu a/b$ on $\mathbb{R}^3\times[0,T/(ab))$, and at corresponding points
\begin{equation}\label{eq:scaling-ratios}
\frac{|v'|^2}{\nu'|\omega'|}=\frac{|v|^2}{\nu|\omega|},\qquad\frac{|v'|^2}{\nu'|\nabla v'|}=\frac{|v|^2}{\nu|\nabla v|}.
\end{equation}
Indeed, the terms $\partial_tv'$, $(v'\cdot\nabla)v'$, $\nabla p'$ and $f'$ all equal $a^2b$ times the corresponding terms of $(v,p,f)$ at $(bx,abt)$, while $\nu'\Delta v'$ equals $\nu'ab^2\Delta v$, and the two agree exactly when $\nu'=\nu a/b$; since $|v'|^2=a^2|v|^2$, $|\omega'|=ab|\omega|$ and $|\nabla v'|=ab|\nabla v|$, the ratios are unchanged.

One should notice also that,a hypothesis written without the viscosity, such as $|v|^2\le|\omega|$, bounds the local Reynolds number $|v|^2/(\nu|\omega|)$ by $1/\nu$ in the units in which it is written. In the rescaled variables it reads $|v'|^2\le(\nu'/\nu)|\omega'|$, which is the same bound: increasing the viscosity by a change of variables strengthens the hypothesis by exactly the same factor. In particular \Cref{thm:main} covers $|v|^2\le|\omega|$ when $\nu>3\sqrt3/4$ in those units, and not otherwise; reaching $C=1$ at $\nu=1$ would require an improvement inside the estimate (\Cref{rem:threshold}), not a change of units.

\section{\texorpdfstring{$\dot H^1$}{H1} criteria}
{For almost every time $t<T^*$, the strong solution, which satisfies $v\in L^2_{\mathrm{loc}}([0,T^*);H^2)$, obeys the energy bound}
\begin{equation}\label{id:EnergyIneq}
\|v(t)\|_{L^2}\le\|v_0\|_{L^2}+T^{1/2}\|f\|_{L^2(0,T;L^2)}.
\end{equation} {and the identity}
\begin{equation}\label{id:H1Energy}
	\frac12\frac{d}{dt}\|\nabla v\|_{L^2}^2+\nu\|\Delta v\|_{L^2}^2=\langle \omega\times v,\Delta v\rangle-\langle f,\Delta v\rangle.
\end{equation}
{Indeed, it suffices to test~\eqref{eq:NS} against $v$ and against $-\Delta v$. The first gives}  
\begin{equation}\frac12\frac{d}{dt}\|v\|_{L^2}^2+\nu\|\nabla v\|_{L^2}^2=\langle f,v\rangle\le\|f\|_{L^2}\|v\|_{L^2},
\end{equation} 
{from which~\eqref{id:EnergyIneq} follows by integration in time. For the second, we write $(v\cdot\nabla)v=\nabla\frac{|v|^2}{2}+\omega\times v$; since $v$ and $\Delta v$ are divergence-free, the pressure gradient and $\nabla\frac{|v|^2}{2}$ drop out, and the nonlinear term reduces to $\langle\omega\times v,\Delta v\rangle$.}

We will also need the following identity, 
\begin{equation}\label{eq:identity}
\int_{\mathbb{R}^3}|\omega|\,v\cdot\Delta v\,dx=-\int_{\mathbb{R}^3}|\omega|^3\,dx-\int_{\mathbb{R}^3}\omega\cdot\big(\nabla|\omega|\times v\big)\,dx .
\end{equation}
which holds for a divergence-free $v\in H^2(\mathbb{R}^3)$ and $\omega=\nabla\times v$.
Indeed, for $\omega\in H^1$ the modulus $|\omega|$ belongs to $H^1$ and satisfies Kato's inequality $|\nabla|\omega||\le|\nabla\omega|$ almost everywhere. In three dimensions $v\in H^2$ gives $v\in L^\infty$ and $\omega,\nabla v\in L^2\cap L^6$, so every integral in~\eqref{eq:identity} converges, and $|\omega|v\in H^1$ because its gradient is bounded pointwise by $|v|\,|\nabla|\omega||+|\omega|\,|\nabla v|\in L^2$. Since $v$ is divergence-free, $\Delta v=-\nabla\times\omega$, and since $\nabla\times v=\omega$, the curl of $|\omega|v$ is $|\omega|\,\omega+\nabla|\omega|\times v$. For fields $A,B\in H^1(\mathbb{R}^3)$ one has $\int A\cdot(\nabla\times B)=\int B\cdot(\nabla\times A)$; applied with $A=|\omega|v$ and $B=\omega$ this gives
\[
\int|\omega|\,v\cdot\Delta v=-\int|\omega|\,v\cdot(\nabla\times\omega)=-\int\omega\cdot\nabla\times\big(|\omega|v\big)=-\int|\omega|^3-\int\omega\cdot\big(\nabla|\omega|\times v\big).
\]
\subsection{An \texorpdfstring{$\dot H^1$}{H1} criterion with an explicit threshold}\label{sec:main}
\begin{theorem}\label{thm:main}
For divergence-free $v_0\in H^1(\mathbb{R}^3)$, $f\in L^2(0,T;L^2(\mathbb{R}^3))$ and $v$ a strong solution of~\eqref{eq:NS} on its maximal interval $[0,T^*)$, $T^*\le T$. under the assumption that there a positive constant $C$ with
\[ C<\frac{4}{3\sqrt3}\]
and a measurable function $M\colon(0,T^*)\to[0,\infty)$ with $\int_0^{T^*}M(t)^{6/5}\,dt<\infty$ such that, for almost every $t\in(0,T^*)$,
\begin{equation}\label{eq:hyp}
|v(x,t)|^2\le C\,\nu\,|\omega(x,t)|\qquad\text{ for almost every } x\in \left\{y,~ |\omega(y,t)|>M(t)\right\}.
\end{equation}
Then $\sup_{0\le t<T^*}\|\nabla v(t)\|_{L^2}<\infty$. Consequently $T^*=T$, and $v$ is a strong solution on $[0,T)$ with $\nabla v\in L^\infty(0,T;L^2)$.
\end{theorem}

\begin{proof}
 At every point $|\omega|v$ is orthogonal to $\omega\times v$, so for every constant $\mu>0$ the non linear term from \eqref{id:H1Energy} equals \[\langle\omega\times v,\Delta v+\frac{\mu}{\nu}|\omega|v\rangle,\] and for every constant $\delta>0$ the Cauchy-Schwarz and Young inequalities give
\[
\langle\omega\times v,\Delta v\rangle\le\frac{\delta\nu}{2}\Big\|\Delta v+\frac{\mu}{\nu}|\omega|v\Big\|_{L^2}^2+\frac{1}{2\delta\nu}\|\omega\times v\|_{L^2}^2 .
\]
Expanding the square and using the identity~\eqref{eq:identity},
\[
\frac{\delta\nu}{2}\Big\|\Delta v+\frac{\mu}{\nu}|\omega|v\Big\|_{L^2}^2=\frac{\delta\nu}{2}\|\Delta v\|_{L^2}^2-\delta\mu\int|\omega|^3-\delta\mu\int\omega\cdot\big(\nabla|\omega|\times v\big)+\frac{\delta\mu^2}{2\nu}\int|\omega|^2|v|^2 .
\]
The cubic term enters with a negative sign and kept as it is. For the cross term, $$\omega\cdot(\nabla|\omega|\times v)=\nabla|\omega|\cdot(v\times\omega),$$ as we have both \(\|\nabla |w|\|\leq|\nabla w|\) and $\|\nabla\omega\|_{L^2}=\|\nabla\times\omega\|_{L^2}=\|\Delta v\|_{L^2}$ we get, for every constant $a>0$, through Cauchy-Schwarz, 
\[
\delta\mu\,\Big|\int\omega\cdot\big(\nabla|\omega|\times v\big)\Big|\le\delta\mu\,\|\Delta v\|_{L^2}\|\omega\times v\|_{L^2}\le\delta\mu\,a\,\nu\,\|\Delta v\|_{L^2}^2+\frac{\delta\mu}{4a\nu}\|\omega\times v\|_{L^2}^2 .
\]
For the forcing term is bound as follows, \[|\langle f,\Delta v\rangle|\le\kappa\nu\|\Delta v\|_{L^2}^2+\|f\|_{L^2}^2/(4\kappa\nu)\] for every constant $\kappa>0$. Together with the inequality $|\omega\times v|\le|\omega||v|$,
\[
\frac{d}{dt}\dfrac{\|\nabla v\|_{L^2}^2}{2}+\nu\Big(1-\frac\delta2-\delta\mu a-\kappa\Big)\|\Delta v\|_{L^2}^2\le-\delta\mu\int|\omega|^3+\frac1\nu\Big(\frac{\delta\mu^2}{2}+\frac{\delta\mu}{4a}+\frac{1}{2\delta}\Big)\int|\omega|^2|v|^2+\frac{\|f\|_{L^2}^2}{4\kappa\nu}.
\]
We now fix $\delta=4/3$ and $\mu={\sqrt3}/{4}$, so that $\delta\mu=1/{\sqrt3}$, ${\delta\mu^2}/{2}=1/8$ and ${1}/{2\delta}=3/8$.\\ 
For $\kappa\in(0,1/3)$ we take $a=(1-3\kappa)/\sqrt3$, which makes the coefficient of $\nu\|\Delta v\|_{L^2}^2$ vanish and gives $${\delta\mu}/{(4a)}={1}/{(4(1-3\kappa))}.$$ 
The bracket in front of $\int|\omega|^2|v|^2$ is then $1/2+{1}/{(4(1-3\kappa))}$, which tends to $3/4$ as $\kappa\to0$. 
Since $C<{4}/{3\sqrt3}=\frac{1}{\sqrt3}\big/\frac34$, we can fix $\kappa\in(0,1/3)$ with
\[
C\Big(\frac12+\frac{1}{4(1-3\kappa)}\Big)\le\frac{1}{\sqrt3},
\]
and obtain
\[
\frac12\frac{d}{dt}\|\nabla v\|_{L^2}^2\le-\frac{1}{\sqrt3}\int|\omega|^3+\frac1\nu\Big(\frac12+\frac{1}{4(1-3\kappa)}\Big)\int|\omega|^2|v|^2+\frac{\|f\|_{L^2}^2}{4\kappa\nu}.
\]
We split the integral of $|\omega|^2|v|^2$ between the sets $\{|\omega|>M(t)\}$ and $\{|\omega|\le M(t)\}$. On the first, hypothesis~\eqref{eq:hyp} gives $|\omega|^2|v|^2\le C\nu|\omega|^3$, so its contribution is at most $$C\Big(\frac12+\frac{1}{4(1-3\kappa)}\Big)\int_{\{|\omega|>M(t)\}}|\omega|^3\le\frac1{\sqrt3}\int|\omega|^3,$$ which the cubic term absorbs. On the second, $|\omega|^2\le M(t)^{6/5}|\omega|^{4/5}$. H\"older's inequality with exponents $5/2$ and $5/3$, the identity $\|\omega\|_{L^2}=\|\nabla v\|_{L^2}$, and the three dimensional Gagliardo-Nirenberg inequality, $$\|v\|_{L^{10/3}}\le C_{\mathrm{GN}}\|v\|_{L^2}^{2/5}\|\nabla v\|_{L^2}^{3/5},$$ gives
\[
\int_{\{|\omega|\le M(t)\}}|\omega|^2|v|^2\le M(t)^{6/5}\|\omega\|_{L^2}^{4/5}\|v\|_{L^{10/3}}^2\le C_{\mathrm{GN}}^2\,M(t)^{6/5}\,\|v\|_{L^2}^{4/5}\,\|\nabla v\|_{L^2}^2 ,
\]
which is linear in $\|\nabla v\|_{L^2}^2$. Hence
\[
\frac12\frac{d}{dt}\|\nabla v(t)\|_{L^2}^2\le\frac{C_{\mathrm{GN}}^2}{\nu}\Big(\frac12+\frac{1}{4(1-3\kappa)}\Big)M(t)^{6/5}\|v(t)\|_{L^2}^{4/5}\|\nabla v(t)\|_{L^2}^2+\frac{\|f(t)\|_{L^2}^2}{4\kappa\nu}.
\]
With the energy estimate \eqref{id:EnergyIneq}, Gronwall's inequality yields, for every $t<T^*$,
\begin{align*}
\|\nabla v(t)\|_{L^2}^2\leq&\Big(\|\nabla v_0\|_{L^2}^2+\frac{1}{2\kappa\nu}\int_0^T\|f(s)\|_{L^2}^2\,ds\Big)\\ &\times\exp\Big(\frac{2C_{\mathrm{GN}}^2}{\nu}\Big(\frac12+\frac{1}{4(1-3\kappa)}\Big)\Big(\|v_0\|_{L^2}+T^{1/2}\|f\|_{L^2(0,T;L^2)}\Big)^{4/5}\int_0^{T^*}M(s)^{6/5}\,ds\Big),
\end{align*}
which is finite and independent of $t$. If $T^*<T$, the local theory for $H^1$ data would force $\|\nabla v(t)\|_{L^2}\to\infty$ as $t\to T^*$, which is excluded; hence $T^*=T$.
\end{proof}

\begin{remark}[The value of the threshold]\label{rem:threshold}
The constant $4/(3\sqrt3)\approx0.770$ is the best that this proof delivers. With a saturated dissipation, that is $\frac\delta2+\delta\mu a=1$, the proof closes for $$C\big(\frac{\delta\mu^2}{2}+\frac{\delta\mu}{4a}+\frac1{2\delta}\big)\le\delta\mu,$$ which, optimizing over $\delta$, $\mu$ and $a$, gives $\delta=\frac43$, $\delta\mu=\frac1{\sqrt3}$, $a=\frac1{\sqrt3}$, where the bracket equals $\frac18+\frac14+\frac38=\frac34$ against a gain of $\frac1{\sqrt3}$.\\ 
Three lossy inequalities had been used: Cauchy-Schwarz, Kato's inequality\footnote{The pointwise inequality \(|\nabla |w||\leq |\nabla w|\).}, and $|\omega\times v|\le|\omega||v|$. A larger constant requires a saving in one of them. Natural sources are the helicity density, since $\|\omega\times v\|_{L^2}^2=\int|\omega|^2|v|^2-\int(v\cdot\omega)^2$; the strictness of Kato's inequality where the direction of the vorticity varies; and the choice of multiplier, since any scalar multiple of $v$ or of $\omega$ may be added to $\Delta v$ without changing $\langle\omega\times v,\Delta v\rangle$.
\end{remark}

\begin{remark}[The role of the threshold $M$]\label{rem:M}
With $M\equiv0$, hypothesis~\eqref{eq:hyp} requires $|v|^2\le C\nu|\omega|$ wherever $\omega\ne0$. By continuity this forces $v$ to vanish on the boundary of the set $\{\omega\ne0\}$, which fails for generic flows: when $\omega(\cdot,t)$ has compact support, $v$ is in general nonzero on the boundary of that support. The threshold removes the difficulty, and only $\int M^{6/5}<\infty$ is used. For instance $M(t)=C_0(T^*-t)^{-\alpha}$ is admissible for every $\alpha<\frac56$. Bounding the low-vorticity part by $M(t)^2\|v\|_{L^2}^2$ instead would require $\int M^2<\infty$, that is $\alpha<\frac12$; we have not found a splitting of this type that requires a smaller power of $M$ than $\frac65$.
\end{remark}
\medskip
\begin{remark}[Physical reading]
Hypothesis~\eqref{eq:hyp} says that the local Reynolds number $|v|^2/(\nu|\omega|)$ stays below $4/(3\sqrt3)$ wherever the vorticity is large. Viscous diffusion then acts on the scale $|v|/|\omega|$ at least as fast as transport, and the theorem states that this suffices to prevent blow-up.
\end{remark}

\subsection{\texorpdfstring{$\dot H^1$}{H1} criteria through the Newton potential of \texorpdfstring{$|\omega|^2$}{|omega|^2}}\label{sec:newton}

Throughout this section $\NP$ denotes the Newton potential on $\mathbb{R}^3$, $${\NP(g)(x)=(-\Delta)^{-1}g(x)}=\frac1{4\pi}\int\frac{g(y)}{|x-y|}\,dy,$$ and $K=3^{-1/2}(2/\pi)^{2/3}<0.43$ is the best constant in the Sobolev inequality $\|h\|_{L^6}\le K\|\nabla h\|_{L^2}$ on $\mathbb{R}^3$~\cite{Tal76}: Talenti's formula $K=(3\pi)^{-1/2}\big(\Gamma(3)/\Gamma(\tfrac32)\big)^{1/3}$ gives this value, and equality holds for $h=(1+|x|^2)^{-1/2}$. In the proof of \Cref{thm:main} a multiple of $|\omega|v$ was added to $\Delta v$ in its pairing with $\omega\times v$. Adding instead a multiple of $\NP(|\omega|^2)v$, the integral of $|\omega|^2|v|^2$ produced by the Cauchy-Schwarz inequality cancels identically, two further terms with a favourable sign appear, and the only term left to control is $\int\NP(|\omega|^2)^2|v|^2$.

\begin{theorem}\label{thm:newton}
Let $v_0\in H^1(\mathbb{R}^3)$ be divergence-free, $f\in L^2(0,T;L^2(\mathbb{R}^3))$, and let $v$ be the strong solution of~\eqref{eq:NS} on its maximal interval $[0,T^*)$, $T^*\le T$. If
\begin{equation}\label{eq:hypnewton}
\int_0^{T^*}\frac{\int_{\mathbb{R}^3}\NP(|\omega|^2)^2\,|v|^2\,dx}{\nu^3\,\|\nabla v\|_{L^2}^2}\,dt<\infty ,
\end{equation}
the integrand being taken equal to $0$ at the times at which $v$ vanishes identically, then $\|\nabla v(t)\|_{L^2}$ stays bounded on $[0,T^*)$, and consequently $T^*=T$. Hypothesis~\eqref{eq:hypnewton} holds in particular when
\begin{equation}\label{eq:hypnewton2}
\int_0^{T^*}\|\omega\|_{L^{12/5}}^2\,\|\omega\|_{L^2}\,\|v\|_{L^3}\,dt<\infty .
\end{equation}
\end{theorem}

\begin{proof}
For $t<T^*$ the strong solution satisfies $v\in L^2_{\mathrm{loc}}([0,T^*);H^2)$. Fix a time $t$ with $v(t)\in H^2$. Then $v\in L^2\cap L^\infty$, $\omega\in H^1\subset L^2\cap L^6$ and $|\omega|^2\in L^1\cap L^3$; consequently $\NP(|\omega|^2)$ is bounded, its gradient belongs to $L^2$, and the integrations by parts below are justified. {We start from the identity~\eqref{id:H1Energy}.} At every point $v$ is orthogonal to $\omega\times v$, so {its nonlinear term satisfies} $\langle\omega\times v,\Delta v\rangle=\big\langle\omega\times v,\Delta v+\nu^{-2}\NP(|\omega|^2)v\big\rangle$; the factor $\nu^{-2}$ makes the two terms of the same dimension. The Cauchy-Schwarz and Young inequalities give
\[
\langle\omega\times v,\Delta v\rangle\le\frac\nu2\Big\|\Delta v+\frac{1}{\nu^2}\NP(|\omega|^2)v\Big\|_{L^2}^2+\frac1{2\nu}\|\omega\times v\|_{L^2}^2 ,
\]
and expanding the square,
\[
\frac\nu2\Big\|\Delta v+\frac{1}{\nu^2}\NP(|\omega|^2)v\Big\|_{L^2}^2=\frac\nu2\|\Delta v\|_{L^2}^2+\frac1\nu\int\NP(|\omega|^2)\,v\cdot\Delta v+\frac1{2\nu^3}\int\NP(|\omega|^2)^2|v|^2 .
\]
The middle term is computed exactly. Since $v\cdot\Delta v=\Delta\frac12|v|^2-|\nabla v|^2$ and $-\Delta\NP(|\omega|^2)=|\omega|^2$, two integrations by parts give
\begin{equation}\label{id:NPgainstdissipation}
\int\NP(|\omega|^2)\,v\cdot\Delta v=-\frac12\int|\omega|^2|v|^2-\int\NP(|\omega|^2)\,|\nabla v|^2 ,
\end{equation}
and both terms on the right are nonpositive. Since $\|\omega\times v\|_{L^2}^2=\int|\omega|^2|v|^2-\int(\omega\cdot v)^2$, the two integrals of $|\omega|^2|v|^2$ cancel exactly, and
\[
\langle\omega\times v,\Delta v\rangle\le\frac\nu2\|\Delta v\|_{L^2}^2-\frac1\nu\int\NP(|\omega|^2)\,|\nabla v|^2-\frac1{2\nu}\int(\omega\cdot v)^2+\frac1{2\nu^3}\int\NP(|\omega|^2)^2|v|^2 .
\]
For the force, \begin{equation}\label{id:fagainstdissipation}
	|\langle f,\Delta v\rangle|\le\frac\nu2\|\Delta v\|_{L^2}^2+\frac1{2\nu}\|f\|_{L^2}^2.
\end{equation} Collecting,
\begin{equation}\label{eq:newton-energy}
\frac12\frac{d}{dt}\|\nabla v\|_{L^2}^2+\frac1\nu\int\NP(|\omega|^2)\,|\nabla v|^2+\frac1{2\nu}\int(\omega\cdot v)^2\le\frac1{2\nu^3}\int\NP(|\omega|^2)^2|v|^2+\frac1{2\nu}\|f\|_{L^2}^2 .
\end{equation}
The two terms on the left after the time derivative are nonnegative and are dropped. Writing $\frac1{2\nu^3}\int\NP(|\omega|^2)^2|v|^2$ as $\frac12\|\nabla v\|_{L^2}^2$ times the integrand of~\eqref{eq:hypnewton}, Gronwall's lemma gives, for every $t<T^*$,
\begin{align*}
\|\nabla v(t)\|_{L^2}^2\le\Big(\|\nabla v_0\|_{L^2}^2+\frac1\nu\int_0^T\|f(s)\|_{L^2}^2\,ds\Big)
\exp\Big(\int_0^{T^*}\frac{\int_{\mathbb{R}^3}\NP(|\omega|^2)^2\,|v|^2\,dx}{\nu^3\,\|\nabla v\|_{L^2}^2}\,ds\Big),
\end{align*}
which is finite and independent of $t$. If $T^*<T$, the local theory for $H^1$ data would force $\|\nabla v(t)\|_{L^2}\to\infty$ as $t\to T^*$, which is excluded; hence $T^*=T$.

It remains to show that~\eqref{eq:hypnewton2} implies~\eqref{eq:hypnewton}. Fix again a time with $v(t)\in H^2$. In $\mathbb{R}^3$ the divergence of $z/|z|$ is $2/|z|$, so $\NP(|\omega|^2)$ is the divergence of the integral of $|\omega(y)|^2(x-y)/(8\pi|x-y|)$, whose modulus is at most $\|\omega\|_{L^2}^2/(8\pi)$ at every point $x$. We write $\int\NP(|\omega|^2)^2|v|^2$ as the integral of this divergence against $\NP(|\omega|^2)|v|^2$. This function is integrable, and so is its gradient $|v|^2\nabla\NP(|\omega|^2)+\NP(|\omega|^2)\nabla|v|^2$, since $|v|^2\in L^1\cap L^2$, $\NP(|\omega|^2)$ is bounded, and $\nabla\NP(|\omega|^2)$ and $\nabla v$ belong to $L^2$. One integration by parts and $|\nabla|v|^2|\le2|v|\,|\nabla v|$ therefore give
\[
\int\NP(|\omega|^2)^2|v|^2\le\frac{\|\omega\|_{L^2}^2}{8\pi}\int\Big(|v|^2\,\big|\nabla\NP(|\omega|^2)\big|+2\,\NP(|\omega|^2)\,|v|\,|\nabla v|\Big)\,dx .
\]
As we have, for $g\in L^{6/5}$, using H\"older then Sobolev embeddings
\[
\|\nabla\NP(g)\|_{L^2}^2=\int g\,\NP(g)\le\|g\|_{L^{6/5}}\|\NP(g)\|_{L^6}\le K\,\|g\|_{L^{6/5}}\,\|\nabla\NP(g)\|_{L^2},
\]
so that $\|\nabla\NP(g)\|_{L^2}\le K\|g\|_{L^{6/5}}$\footnote{ The best constant is given by duality with the Sobolev embedding.}, and $\|\NP(g)\|_{L^6}\le K^2\|g\|_{L^{6/5}}$, we get for the first integral, with $g=|\omega|^2$ then  Cauchy-Schwarz, that it is smaller than    $$K\|v\|_{L^4}^2\|\omega\|_{L^{12/5}}^2,$$ using H\"older's  with exponents $6$, $3$ and $2$ bounds the second one by $$2K^2\|\omega\|_{L^{12/5}}^2\|v\|_{L^3}\|\nabla v\|_{L^2}.$$ Finally $\|v\|_{L^4}^2\le\|v\|_{L^3}\|v\|_{L^6}\le K\|v\|_{L^3}\|\nabla v\|_{L^2}$ and $\|\nabla v\|_{L^2}=\|\omega\|_{L^2}$, so that
\[
\int\NP(|\omega|^2)^2|v|^2\le\frac{3K^2}{8\pi}\,\|\omega\|_{L^2}^2\,\|\omega\|_{L^{12/5}}^2\,\|\omega\|_{L^2}\,\|v\|_{L^3}.
\]
Dividing by $\nu^3\|\nabla v\|_{L^2}^2=\nu^3\|\omega\|_{L^2}^2$, the integrand of~\eqref{eq:hypnewton} is at most ${3K^2}/{8\pi\nu^3}$ times the integrand of~\eqref{eq:hypnewton2}, which proves the last assertion.
\end{proof}
\medskip
\begin{remark}[Scaling and the reach of the norm condition]\label{rem:newton-scaling}
Under $v\mapsto\lambda v(\lambda x,\lambda^2t)$ the integrands of~\eqref{eq:hypnewton} and~\eqref{eq:hypnewton2} are multiplied by $\lambda^2$ and $dt$ by $\lambda^{-2}$, so both hypotheses are critical. Condition~\eqref{eq:hypnewton2} is only a sufficient condition for~\eqref{eq:hypnewton}, and its factors are close to known criteria: a bound on $\sup_t\|v\|_{L^3}$ alone already gives regularity in the unforced case~\cite{ESS03}, so~\eqref{eq:hypnewton2} may adds something only when $\|v\|_{L^3}$ is not known to be bounded, as it is along the forced solution (\Cref{rem:newton-forced}). This is unavoidable for any bound of~\eqref{eq:hypnewton} by norms. The integrations by parts bring out $\sup_t\||v|^2\|_{L^{6/5}}=\sup_t\|v\|_{L^{12/5}}^2$ and $\int\||\omega|^2\|_{L^{6/5}}\,dt=\int\|\omega\|_{L^{12/5}}^2\,dt$, which are multiplied by $\lambda^{-1/2}$ under the scaling, halfway between the energy level and the critical level; a critical bound must pair them with factors multiplied by $\lambda^{1/2}$, and such factors are themselves close to regularity criteria. The content of \Cref{thm:newton} lies in~\eqref{eq:hypnewton} and in the exact cancellation of $\int|\omega|^2|v|^2$ in its proof. Without a force the whole dissipation is available: with the multiplier $\frac1{4\nu^2}\NP(|\omega|^2)v$ and the Young inequality $\langle\omega\times v,\Delta v\rangle\le\nu\|\Delta v+\frac1{4\nu^2}\NP(|\omega|^2)v\|_{L^2}^2+\frac1{4\nu}\|\omega\times v\|_{L^2}^2$, the same computation gives $\frac{d}{dt}\|\nabla v\|_{L^2}^2\le\frac1{8\nu^3}\int\NP(|\omega|^2)^2|v|^2$, so the exponent in the Gronwall bound is divided by $8$.
\end{remark}
\medskip
\begin{remark}[The forced blow-up]\label{rem:newton-forced}
With the scalings recalled in \Cref{sec:construction} of the present paper, at leading order and up to the logarithmic factors of the pulse annulus, the forced solution of~\cite{OpenAI26} violates~\eqref{eq:hypnewton2}, and only by the anomalous factor. At height $z$ the radial integral of $|v|^3$ is of order $q^{-1/2-3h}$, where $q$ is the concentration scale, comparable to $(1-t)+|z|^{2/(1-2h)}$; since $1+6h>1-2h$, its integral in $z$ is dominated by the core, of volume comparable to $(1-t)^{3/2-h}$, and the same holds for $|\omega|^{12/5}$ and $|\omega|^2$. This gives
	\[
	\int|v|^3\asymp(1-t)^{-4h},\qquad\int|\omega|^{12/5}\asymp(1-t)^{-\frac9{10}-\frac{17h}5},\qquad\int|\omega|^2\asymp(1-t)^{-\frac12-3h},
	\]
	so that $\|\omega\|_{L^{12/5}}^2\|\omega\|_{L^2}\|v\|_{L^3}$ is of order $(1-t)^{-1-17h/3}$, which is not integrable on $(0,1)$. Here $\|v\|_{L^3}$ is unbounded, but only like $(1-t)^{-4h/3}$: this is the situation in which~\eqref{eq:hypnewton2} goes beyond a bound on $\sup_t\|v\|_{L^3}$. The same solution keeps the factors in $L^{12/5}$ under control. For $|v|^{12/5}$ the radial integral is of order $q^{-1/5-12h/5}$, whose integral in $z$ converges even at $t=1$ because $2(1+12h)<5(1-2h)$, so $\sup_{t<1}\|v(t)\|_{L^{12/5}}<\infty$; and $\int_0^1\|\omega\|_{L^{12/5}}^2\,dt<\infty$ because $\frac34+\frac{17h}6<1$. These two quantities alone therefore do not imply regularity for the forced equations.
\end{remark}

The cancellation can also be produced on the other side of the pairing. Since $\Delta v$ is divergence-free, it is orthogonal to every gradient in $L^2$, so a gradient may be added to $\omega\times v$ without changing $\langle\omega\times v,\Delta v\rangle$. The gradient $\nabla\NP(|\omega|^2)$ has the dimension of $\omega\times v$; with it the cross term again cancels $\int|\omega|^2|v|^2$ exactly, and what remains involves the vorticity alone.

\begin{theorem}\label{thm:coulomb}
Let $v_0\in H^1(\mathbb{R}^3)$ be divergence-free, $f\in L^2(0,T;L^2(\mathbb{R}^3))$, and let $v$ be the strong solution of~\eqref{eq:NS} on its maximal interval $[0,T^*)$, $T^*\le T$. If
\begin{equation}\label{eq:hypcoulomb}
\int_0^{T^*}\frac{\int_{\mathbb{R}^3}\NP(|\omega|^2)\,|\omega|^2\,dx}{\nu\,\|\nabla v\|_{L^2}^2}\,dt<\infty ,
\end{equation}
the integrand being taken equal to $0$ at the times at which $v$ vanishes identically, then $\|\nabla v(t)\|_{L^2}$ stays bounded on $[0,T^*)$, and consequently $T^*=T$. Hypothesis~\eqref{eq:hypcoulomb} holds in particular when
\[
\int_0^{T^*}\frac{\|\omega\|_{L^{12/5}}^4}{\|\omega\|_{L^2}^2}\,dt<\infty ,
\]
and therefore when $\int_0^{T^*}\|\omega\|_{L^3}^2\,dt<\infty$.
\end{theorem}

\begin{proof}
Fix a time $t$ with $v(t)\in H^2$. As in the proof of \Cref{thm:newton}, $|\omega|^2\in L^1\cap L^3$, $\NP(|\omega|^2)$ is bounded and its gradient belongs to $L^2$, which justifies the integrations by parts below. {We start again from the identity~\eqref{id:H1Energy}.} Since $\Delta v$ is divergence-free and $\nabla\NP(|\omega|^2)\in L^2$, the two are orthogonal, so $\langle\omega\times v,\Delta v\rangle=\langle\omega\times v+\nabla\NP(|\omega|^2),\Delta v\rangle$; the two terms of the sum have the same dimension, and no power of $\nu$ is needed. The Cauchy-Schwarz and Young inequalities give
\[
\langle\omega\times v,\Delta v\rangle\le\frac\nu2\|\Delta v\|_{L^2}^2+\frac1{2\nu}\big\|\omega\times v+\nabla\NP(|\omega|^2)\big\|_{L^2}^2 ,
\]
and the square is computed exactly. Its first term is $\|\omega\times v\|_{L^2}^2=\int|\omega|^2|v|^2-\int(\omega\cdot v)^2$, and its last term is $\|\nabla\NP(|\omega|^2)\|_{L^2}^2=\int\NP(|\omega|^2)\,|\omega|^2$. For the cross term, one integration by parts and $\nabla\cdot(\omega\times v)=v\cdot(\nabla\times\omega)-|\omega|^2=-v\cdot\Delta v-|\omega|^2$ give
\[
2\big\langle\omega\times v,\nabla\NP(|\omega|^2)\big\rangle=2\int\NP(|\omega|^2)\,v\cdot\Delta v+2\int\NP(|\omega|^2)\,|\omega|^2.
\]
With the identity \eqref{id:NPgainstdissipation} for $\int\NP(|\omega|^2)\,v\cdot\Delta v$, the two integrals of $|\omega|^2|v|^2$ cancel exactly, and
\[
\big\|\omega\times v+\nabla\NP(|\omega|^2)\big\|_{L^2}^2=3\int\NP(|\omega|^2)\,|\omega|^2-2\int\NP(|\omega|^2)\,|\nabla v|^2-\int(\omega\cdot v)^2 .
\]
 Collect these bounds together with \eqref{id:fagainstdissipation} for the force, and write $\frac1\nu\int\NP(|\omega|^2)|\nabla v|^2$ as $$\frac1\nu\int\NP(|\omega|^2)\big(|\nabla v|^2-\frac12|\omega|^2\big)+\frac1{2\nu}\int\NP(|\omega|^2)|\omega|^2,$$ you get
\begin{equation}\label{eq:coulomb-energy}
\frac12\frac{d}{dt}\|\nabla v\|_{L^2}^2+\frac1\nu\int\NP(|\omega|^2)\Big(|\nabla v|^2-\frac12|\omega|^2\Big)+\frac1{2\nu}\int(\omega\cdot v)^2\le\frac1\nu\int\NP(|\omega|^2)\,|\omega|^2+\frac1{2\nu}\|f\|_{L^2}^2 .
\end{equation}
By~\eqref{eq:split} and the positivity of $\NP(|\omega|^2)$, the two terms on the left after the time derivative are nonnegative, and they are dropped. Writing $\frac1\nu\int\NP(|\omega|^2)|\omega|^2$ as $\|\nabla v\|_{L^2}^2$ times the integrand of~\eqref{eq:hypcoulomb}, Gronwall's lemma gives, for every $t<T^*$,
\begin{align*}
\|\nabla v(t)\|_{L^2}^2\le\Big(\|\nabla v_0\|_{L^2}^2+\frac1\nu\int_0^T\|f(s)\|_{L^2}^2\,ds\Big)
\exp\Big(2\int_0^{T^*}\frac{\int_{\mathbb{R}^3}\NP(|\omega|^2)\,|\omega|^2\,dx}{\nu\,\|\nabla v\|_{L^2}^2}\,ds\Big),
\end{align*}
which is finite and independent of $t$, and $T^*=T$ as in the proof of \Cref{thm:newton}. For the last assertion, as in that proof, $\int\NP(|\omega|^2)|\omega|^2=\|\nabla\NP(|\omega|^2)\|_{L^2}^2\le K^2\|\omega\|_{L^{12/5}}^4$, and H\"older's inequality with $\frac5{12}=\frac12\cdot\frac12+\frac12\cdot\frac13$ gives $\|\omega\|_{L^{12/5}}^2\le\|\omega\|_{L^2}\|\omega\|_{L^3}$. Since $\|\nabla v\|_{L^2}=\|\omega\|_{L^2}$, the integrand of~\eqref{eq:hypcoulomb} is at most $\frac{K^2}{\nu}\|\omega\|_{L^{12/5}}^4\big/\|\omega\|_{L^2}^2$, which is at most $\frac{K^2}{\nu}\|\omega\|_{L^3}^2$.
\end{proof}

\begin{remark}[Scaling and comparison]\label{rem:coulomb}
Under $v\mapsto\lambda v(\lambda x,\lambda^2t)$, $\int\NP(|\omega|^2)|\omega|^2$ is multiplied by $\lambda^3$ and $\|\nabla v\|_{L^2}^2$ by $\lambda$, so the integrand of~\eqref{eq:hypcoulomb} is multiplied by $\lambda^2$ and the hypothesis is critical, without smallness. Its numerator is the Newtonian self-energy of $|\omega|^2$,
\[
\int_{\mathbb{R}^3}\NP(|\omega|^2)\,|\omega|^2\,dx=\frac1{4\pi}\iint\frac{|\omega(x)|^2\,|\omega(y)|^2}{|x-y|}\,dx\,dy ,
\]
so the hypothesis involves the vorticity alone. By the last assertion of \Cref{thm:coulomb} it is implied by $\int\|\omega\|_{L^3}^2\,dt<\infty$, the case of spatial exponent $3$ of the critical vorticity condition of~\cite{BdV95}. The intermediate condition $\int\|\omega\|_{L^{12/5}}^4\|\omega\|_{L^2}^{-2}\,dt<\infty$ does not compare in an obvious way with the case of spatial exponent $\frac{12}{5}$, $\int\|\omega\|_{L^{12/5}}^{8/3}\,dt<\infty$.\\ Hypotheses~\eqref{eq:hypnewton} and~\eqref{eq:hypcoulomb} are not comparable either. At each point the density $\NP(|\omega|^2)^2|v|^2/\nu^3$ of the first is the density $\NP(|\omega|^2)|\omega|^2/\nu$ of the second times $\NP(|\omega|^2)|v|^2/(\nu^2|\omega|^2)$, the square of a local Reynolds number built on the length $\NP(|\omega|^2)^{1/2}/|\omega|$; \Cref{thm:newton} is the sharper of the two where this number is small, and \Cref{thm:coulomb} where it is large. Two variants follow from~\eqref{eq:coulomb-energy}. Since $|\nabla v|^2\ge\frac12|\omega|^2$, the condition with $|\nabla v|^2$ in place of $|\omega|^2$ in the numerator of~\eqref{eq:hypcoulomb} also suffices. And since the second term on the left of~\eqref{eq:coulomb-energy} was dropped, only the excess of rotation over strain counts: by~\eqref{eq:split}, the numerator of~\eqref{eq:hypcoulomb} may be replaced by the positive part of $\int\NP(|\omega|^2)\big(\frac32|\omega|^2-|\nabla v|^2\big)$.
\end{remark}
\medskip
\begin{remark}[The forced blow-up]\label{rem:coulomb-forced}
With the scalings recalled in \Cref{sec:construction} of the present paper, at leading order, the forced solution of~\cite{OpenAI26} violates~\eqref{eq:hypcoulomb}, and only by the anomalous factor. In the core, of radius comparable to $(1-t)^{1/2}$ and length comparable to $(1-t)^{1/2-h}$, $|\omega|^2$ is of order $(1-t)^{-2-2h}$. The Newton potential of a uniform density on a slender column is of order the density times the square of the radius, times one plus the logarithm of the ratio of length to radius, here $h\ln\frac1{1-t}$; so $\NP(|\omega|^2)$ is of order $(1-t)^{-1-2h}\big(1+h\ln\frac1{1-t}\big)$ in the core. The self-energy is dominated by the core, as are the integrals of \Cref{rem:newton-forced}, and
\[
\int\NP(|\omega|^2)\,|\omega|^2\asymp(1-t)^{-\frac32-5h}\Big(1+h\ln\frac1{1-t}\Big),\qquad\int|\omega|^2\asymp(1-t)^{-\frac12-3h},
\]
so that the integrand of~\eqref{eq:hypcoulomb} is of order $(1-t)^{-1-2h}\big(1+h\ln\frac1{1-t}\big)$, which is not integrable on $(0,1)$. For comparison, $\|\omega\|_{L^3}^2\asymp(1-t)^{-1-8h/3}$ by \Cref{sec:forced}, in agreement with the last assertion of \Cref{thm:coulomb}.
\end{remark}

\section{A maximum principle for \texorpdfstring{$\frac12|v|^2+\lambda|\omega|$}{1/2|v|^2+lambda|omega|}}\label{sec:max}

Throughout this section $v_0\in H^3(\mathbb{R}^3)$ is divergence-free, $f\in C^\infty_c(\mathbb{R}^3\times[0,T])$, and $v$ is the smooth solution of~\eqref{eq:NS} on its maximal interval $[0,T^*)$. For $t<T^*$, $v(t)\in H^3$, so $v$, $\nabla v$ and $\omega$ are continuous and tend to zero at infinity. The pressure is normalised by $p=\Phi*\big(\operatorname{tr}((\nabla v)^2)-\nabla\cdot f\big)$, where $\Phi(z)=1/(4\pi|z|)$ is the fundamental solution of $-\Delta$; this is the representation obtained by taking the divergence of~\eqref{eq:NS}, since $\nabla\cdot((v\cdot\nabla)v)=\operatorname{tr}((\nabla v)^2)$ for divergence-free $v$. For a constant $\lambda>0$ we write, following the quantity used by the author,
\[
m=\frac12|v|^2+\lambda|\omega| .
\]

\subsection{The pointwise inequality}

\begin{proposition}\label{prop:pointwise}
For all constants $\lambda,\varepsilon>0$, on $\mathbb{R}^3\times(0,T^*)$,
\begin{equation}\label{eq:pointwise}
(\partial_t+v\cdot\nabla-\nu\Delta)\Big(\frac12|v|^2+\lambda\big(|\omega|^2+\varepsilon^2\big)^{1/2}\Big)\le-\Big(\nu-\frac{\lambda}{\sqrt3}\Big)|\nabla v|^2-v\cdot\nabla p+f\cdot v+\lambda|\nabla\times f| .
\end{equation}
\end{proposition}

\begin{proof}
Taking the scalar product of~\eqref{eq:NS} with $v$ and using $v\cdot\Delta v=\Delta\frac12|v|^2-|\nabla v|^2$ gives $$(\partial_t+v\cdot\nabla-\nu\Delta)\frac12|v|^2=-\nu|\nabla v|^2-v\cdot\nabla p+f\cdot v.$$ The vorticity satisfies $\partial_t\omega+(v\cdot\nabla)\omega-\nu\Delta\omega=(\omega\cdot\nabla)v+\nabla\times f$. For $(|\omega|^2+\varepsilon^2)^{1/2}$, the chain rule gives
\[
(\partial_t+v\cdot\nabla-\nu\Delta)\big(|\omega|^2+\varepsilon^2\big)^{1/2}=\frac{\omega\cdot(\omega\cdot\nabla)v+\omega\cdot\nabla\times f}{(|\omega|^2+\varepsilon^2)^{1/2}}-\nu\,\frac{|\nabla\omega|^2-\big|\nabla\frac12|\omega|^2\big|^2\big/(|\omega|^2+\varepsilon^2)}{(|\omega|^2+\varepsilon^2)^{1/2}} .
\]
The last term is nonpositive, since $|\nabla\frac12|\omega|^2|\le|\omega||\nabla\omega|\le(|\omega|^2+\varepsilon^2)^{1/2}|\nabla\omega|$. In the first, $\omega\cdot(\omega\cdot\nabla)v=\omega\cdot S\omega$\footnote{Here \(S\) is the strain tensor as the anti symmetric part yields $w\times w=0$.} is at most the largest eigenvalue of $S$ times $|\omega|^2$, and this eigenvalue is nonnegative because $S$ is trace-free. For a trace-free symmetric matrix the largest eigenvalue is minus the sum of the other two, whose squares add up to at least half its square, so the squared norm of the matrix is at least $\frac32$ times its square: the largest eigenvalue is at most $\sqrt{2/3}\,|S|$. By~\eqref{eq:split}, $|S||\omega|=\sqrt2|S|\cdot|\omega|/\sqrt2\le\big(|S|^2+\frac12|\omega|^2\big)/\sqrt2=|\nabla v|^2/\sqrt2$. Hence the stretching contribution is at most $\sqrt{2/3}\,|S||\omega|\le|\nabla v|^2/\sqrt3$, the forcing contribution is at most $|\nabla\times f|$, and adding $\lambda$ times this to the energy part gives~\eqref{eq:pointwise}.
\end{proof}

Vortex stretching is thus absorbed pointwise by viscous dissipation for every $\lambda<\sqrt3\,\nu$. The only term without a sign is the pressure work.

\subsection{The pressure at a point}

\begin{lemma}\label{lem:weighted}
There is an absolute constant $C_3$ such that, for every divergence-free $v$ with $\nabla v\in L^2(\mathbb{R}^3)$ and $\omega\in L^2\cap L^\infty(\mathbb{R}^3)$, every $x\in\mathbb{R}^3$ and every $r_0>0$,
\begin{equation}\label{eq:weighted}
\int_{B(x,r_0)}\frac{|\nabla v(y)|^2}{|x-y|^2}\,dy\le C_3\,r_0\,\|\omega\|_{L^\infty}^2+\frac{64\pi}{3}\,r_0\,\Big|\frac{1}{|B(x,r_0)|}\int_{B(x,r_0)}\nabla v(z)\,dz\Big|^2 .
\end{equation}
If moreover $v\in L^\infty$, the last term is at most $192\pi\,\|v\|_{L^\infty}^2/r_0$.
\end{lemma}

\begin{proof}
Write $\NP$ for the inverse of $-\Delta$ on $\mathbb{R}^3$. Since $-\Delta v=\nabla\times\omega$, we have $v=\nabla\times\NP(\omega)$. Fix a ball $B_r(x_0)$ and split $\omega=\omega\mathbf 1_{B_{2r}(x_0)}+\omega\mathbf 1_{\mathbb{R}^3\setminus B_{2r}(x_0)}$. 
For any $U\in L^2(\mathbb{R}^3;\mathbb{R}^3)$, the Fourier multiplier of $\nabla\nabla\times\NP$ has squared norm $|\xi\times\hat U|^2/|\xi|^2\le|\hat U|^2$, so the local piece satisfies $$\int|\nabla\nabla\times\NP(\omega\mathbf 1_{B_{2r}(x_0)})|^2\le\int_{B_{2r}(x_0)}|\omega|^2\le\frac{32\pi}{3}r^3\|\omega\|_{L^\infty}^2.$$ 
For $y\in B_r(x_0)$ the far piece is given by an absolutely convergent integral of $\omega(z)$ against the second derivatives of $1/|y-z|$, over $|z-x_0|\ge2r$. The third derivatives of $1/|y-z|$ are bounded by $c|y-z|^{-4}$, and every point of the segment from $x_0$ to $y$ lies at distance at least $\frac12|x_0-z|$ from $z$; by the mean value theorem the far piece at $y$ differs from its value at $x_0$ by at most $$\frac{1}{4\pi}\,16cr\,\|\omega\|_{L^\infty}\int_{|z-x_0|\ge2r}|x_0-z|^{-4}\,dz=8c\,\|\omega\|_{L^\infty}=:C_0\|\omega\|_{L^\infty}.$$ The factor $r$ cancels exactly, which is why no logarithm appears. Since the ball average of $\nabla v$ minimises $\int_{B_r(x_0)}|\nabla v-A|^2$ over constant matrices $A$, taking for $A$ the far piece at $x_0$ gives
\[
\int_{B_r(x_0)}\Big|\nabla v-\frac{1}{|B_r(x_0)|}\int_{B_r(x_0)}\nabla v\Big|^2\le C_1\,r^3\,\|\omega\|_{L^\infty}^2,\qquad C_1=\frac{64\pi}{3}+\frac{8\pi}{3}C_0^2,
\]
for every $x_0$ and $r$: the $L^2$ form of the bound of $\nabla v$ in $BMO$ by $\|\omega\|_{L^\infty}$. By Cauchy-Schwarz and this bound at radius $2\rho$, the averages of $\nabla v$ over $B_\rho(x)$ and $B_{2\rho}(x)$ differ by at most $C_2\|\omega\|_{L^\infty}$ with $C_2=(6C_1/\pi)^{1/2}$, so by telescoping, the average over $B_{2^{-k}r_0}(x)$ has norm at most the norm of the average over $B_{r_0}(x)$ plus $kC_2\|\omega\|_{L^\infty}$. The exact identity $$\int_B|\nabla v|^2=\int_B|\nabla v-\frac1{|B|}\int_B\nabla v|^2+|B|\,|\frac1{|B|}\int_B\nabla v|^2$$ then gives
\[
\int_{B_{2^{-k}r_0}(x)}|\nabla v|^2\le2^{-3k}r_0^3\Big[C_1\|\omega\|_{L^\infty}^2+\frac{8\pi}{3}\Big|\frac{1}{|B(x,r_0)|}\int_{B(x,r_0)}\nabla v\Big|^2+\frac{8\pi}{3}k^2C_2^2\|\omega\|_{L^\infty}^2\Big].
\]
On the shell $2^{-k-1}r_0\le|x-y|<2^{-k}r_0$ the weight is at most $4^{k+1}r_0^{-2}$. Summing over $k\ge0$ with $\sum_k4\cdot2^{-k}=8$ and $\sum_k4k^2\,2^{-k}=24$ gives~\eqref{eq:weighted} with $C_3=8C_1+64\pi C_2^2$. Finally, by the divergence theorem the average of $\nabla v$ over $B(x,r_0)$ equals $$\frac{1}{|B(x,r_0)|}\int_{\partial B(x,r_0)}v\otimes n\,dS,$$ whose norm is at most ${3}\|v\|_{L^\infty}/r_0$, and ${64\pi}\cdot9/{3}=192\pi$.
\end{proof}

\begin{remark}
Running the same argument from a top scale $R$ and bounding the average there by $|B_R|^{-1/2}\|\nabla v\|_{L^2}=|B_R|^{-1/2}\|\omega\|_{L^2}$ gives the global bound $$\sup_x\int_{\mathbb{R}^3}|\nabla v(y)|^2|x-y|^{-2}\,dy\le(17+C_3)\|\omega\|_{L^\infty}^{4/3}\|\omega\|_{L^2}^{2/3},$$ with $R=(\|\omega\|_{L^2}/\|\omega\|_{L^\infty})^{2/3}$. We do not use it here.
\end{remark}

\begin{corollary}\label{cor:pressure}
There is an absolute constant $c_0$ such that, for every $x\in\mathbb{R}^3$ and $t\in(0,T^*)$,
\[
|\nabla p(x,t)|\le c_0\,\|v(t)\|_{L^\infty}\|\omega(t)\|_{L^\infty}+A_f,\qquad A_f=\sup_{0\le s\le T}\big\|\nabla\Phi*\nabla\cdot f(s)\big\|_{L^\infty}<\infty .
\]
\end{corollary}

\begin{proof}\sloppy
The force contribution is at most $A_f$. 
If $\omega(t)\equiv0$ then $v(t)\equiv0$; otherwise set $r_0=\|v\|_{L^\infty}/\|\omega\|_{L^\infty}$, and let $\chi$ be a smooth radial cutoff equal to $1$ on $B(0,r_0)$, vanishing outside $B(0,2r_0)$, with $|\nabla\chi|\le2/r_0$ and $|\nabla^2\chi|\le c/r_0^2$. Since $|\operatorname{tr}((\nabla v)^2)|\le|\nabla v|^2$ and $|\nabla\Phi(z)|=1/(4\pi|z|^2)$, the near part of $\nabla\Phi*\operatorname{tr}((\nabla v)^2)$ at $x$ is bounded, by \Cref{lem:weighted} at radius $2r_0$, by
\[
\frac1{4\pi}\int_{B(x,2r_0)}\frac{|\nabla v(y)|^2}{|x-y|^2}\,dy\le\frac{C_3}{2\pi}\,r_0\|\omega\|_{L^\infty}^2+\frac{24}{r_0}\|v\|_{L^\infty}^2 .
\]
For the far part we write $\operatorname{tr}((\nabla v)^2)=\nabla\cdot\nabla\cdot(v\otimes v)$ and integrate by parts twice against the smooth kernel $(1-\chi(x-\cdot))\nabla\Phi(x-\cdot)$. Its second derivatives are bounded by $c_K|x-y|^{-4}$ and vanish for $|x-y|<r_0$, with $c_K$ independent of $r_0$ by scaling, so the far part is at most $$c_K\|v\|_{L^\infty}^2\int_{|z|\ge r_0}|z|^{-4}\,dz=4\pi c_K\|v\|_{L^\infty}^2/r_0.$$ With $r_0=\|v\|_{L^\infty}/\|\omega\|_{L^\infty}$ the three bounds add up to $c_0\|v\|_{L^\infty}\|\omega\|_{L^\infty}$ with $$c_0=24+4\pi c_K+{C_3}/{2\pi}.$$
\end{proof}
\subsection{The maximum principle}

\begin{theorem}\label{thm:max}
Let $c_0$ be the constant of \Cref{cor:pressure}, and set
\[
\lambda=\frac{\nu}{32c_0+1},\qquad C_*=\Big(\frac{1-1/(\sqrt3(32c_0+1))}{2c_0(32c_0+1)}\Big)^{1/2}.
\]
Assume that there are constants $\eta\in(0,1)$ and $C\in[0,C_*)$ such that, for every $t\in(0,T^*)$,
\begin{equation}\label{eq:hypmax}
|v(x,t)|^2\le C\,\nu\,|\nabla v(x,t)|
\end{equation}
at every point $x$ where both $$m(x,t)\ge(1-\eta)\sup_ym(y,t)$$ and $$|v(x,t)|^2\ge(1-\eta)\sup_ym(y,t).$$ Then $\sup_{t<T^*}\sup_x m(x,t)<\infty$. Consequently $\sup_{t<T^*}\|v(t)\|_{L^\infty}<\infty$ and $T^*=T$.
\end{theorem}

\begin{proof}\sloppy
Fix $\varepsilon>0$. For each $t<T^*$ the function $\frac12|v|^2+\lambda(|\omega|^2+\varepsilon^2)^{1/2}-\lambda\varepsilon$ is continuous, nonnegative and tends to zero at infinity, so the maximum
\[
\Mx
\]
is attained. It is locally Lipschitz in $t$, and by the standard argument for maxima of smooth families its upper right Dini derivative is at most the time derivative of $\frac12|v|^2+\lambda(|\omega|^2+\varepsilon^2)^{1/2}$ at some maximum point $x$. At such a point the gradient vanishes and the Laplacian is nonpositive, so~\eqref{eq:pointwise} gives
\begin{multline*}
D^+\Mx\\
\le-\Big(\nu-\frac{\lambda}{\sqrt3}\Big)|\nabla v(x,t)|^2-v(x,t)\cdot\nabla p(x,t)+|f(x,t)||v(x,t)|+\lambda|\nabla\times f(x,t)| .
\end{multline*}
Everywhere, $\frac12|v|^2$ and $\lambda|\omega|$ are at most this maximum, so $\|v\|_{L^\infty}^2$ is at most twice it and $\|\omega\|_{L^\infty}$ at most its quotient by $\lambda$. \Cref{cor:pressure} then gives
\begin{multline*}
|v(x,t)\cdot\nabla p(x,t)|\le\frac{2c_0}{\lambda}\Big(\Mx\Big)^2\\
+A_f\Big(2\Mx\Big)^{1/2}.
\end{multline*}
Assume the maximum is at least $\lambda\varepsilon\max\{8,\,1+1/\eta\}$. Either $|v(x,t)|^2$ is smaller than the maximum or it is not.

In the first case $\lambda(|\omega(x,t)|^2+\varepsilon^2)^{1/2}$ exceeds half the maximum, so $\lambda|\omega(x,t)|$ is at least half the maximum minus $\lambda\varepsilon$, hence at least $\frac38$ of the maximum, and by~\eqref{eq:split}
\[
|\nabla v(x,t)|^2\ge\frac12|\omega(x,t)|^2\ge\frac{9}{128\lambda^2}\Big(\Mx\Big)^2 .
\]
With $\nu/\lambda=32c_0+1$,
\begin{align*}
-\Big(\nu-\frac{\lambda}{\sqrt3}\Big)|\nabla v(x,t)|^2&+\frac{2c_0}{\lambda}\Big(\Mx\Big)^2\\
&\le-\frac1\lambda\Big(\frac{9(32c_0+1-1/\sqrt3)}{128}-2c_0\Big)\Big(\Mx\Big)^2\\
&\le-\frac{c_0}{4\lambda}\Big(\Mx\Big)^2 .
\end{align*}
No hypothesis is used in this case.

In the second case $|v(x,t)|^2$ is at least the maximum, which is at least $\sup_ym(y,t)$. Moreover $m(x,t)$ is at least the maximum minus $\lambda\varepsilon$, hence at least $\sup_ym(y,t)-\lambda\varepsilon$, which is at least $(1-\eta)\sup_ym(y,t)$ because $\sup_ym(y,t)$ is at least the maximum minus $\lambda\varepsilon$, hence at least $\lambda\varepsilon/\eta$. Hypothesis~\eqref{eq:hypmax} therefore applies at $x$, and
\[
|\nabla v(x,t)|^2\ge\frac{|v(x,t)|^4}{C^2\nu^2}\ge\frac{1}{C^2\nu^2}\Big(\Mx\Big)^2 .
\]
Hence
\begin{align*}
-\Big(\nu-\frac{\lambda}{\sqrt3}\Big)|\nabla v(x,t)|^2&+\frac{2c_0}{\lambda}\Big(\Mx\Big)^2\\
&\le-\frac1\nu\Big(\frac{1-1/(\sqrt3(32c_0+1))}{C^2}-2c_0(32c_0+1)\Big)\\
&\qquad\times\Big(\Mx\Big)^2,
\end{align*}
and the bracket is positive because $C<C_*$.

In both cases, with $k>0$ the smaller of $\frac{c_0}{4}(32c_0+1)$ and the last bracket,
\begin{multline*}
D^+\Mx\le-\frac{k}{\nu}\Big(\Mx\Big)^2\\
+\big(A_f+\|f\|_{L^\infty}\big)\Big(2\Mx\Big)^{1/2}+\lambda\|\nabla\times f\|_{L^\infty} .
\end{multline*}
The right-hand side is negative as soon as the maximum exceeds a value $K_f$ that depends only on $k$, $\nu$, $\lambda$ and the force. Hence, for every $t<T^*$, the maximum is at most the largest of its value at $t=0$, of $\lambda\varepsilon\max\{8,1+1/\eta\}$ and of $K_f$. Letting $\varepsilon\to0$, $\sup_xm(x,t)$ is bounded on $[0,T^*)$ by the largest of $\sup_xm(x,0)$ and $K_f$. Then $\|v(t)\|_{L^\infty}^2\le2\sup_xm(x,t)$ is bounded on $[0,T^*)$, and a velocity bounded in $L^\infty$ cannot blow up, so $T^*=T$.
\end{proof}

\begin{remark}[On the hypothesis and the constant]\label{rem:max}
Hypothesis~\eqref{eq:hypmax} bounds the local Reynolds number $|v|^2/(\nu|\nabla v|)$, and only at the near-maximum points of $m$ where the velocity part dominates. At the near-maximum points where the vorticity part dominates nothing is assumed: viscosity alone defeats both stretching and pressure there. By~\eqref{eq:scaling-ratios} the threshold $C_*\approx1/(8c_0)$ cannot be changed by a change of variables; in units in which the hypothesis is written as $|v|^2\le|\nabla v|$, \Cref{thm:max} applies when $\nu>1/C_*$. The constants have not been optimized. Bounding the third derivatives of $1/|z|$ by $\sqrt{90}\,|z|^{-4}$ and taking a radial cutoff with $|\eta'|\le2$, $|\eta''|\le8$, one finds $C_3\approx1.89\times10^{7}$ and $c_K\le6.1$, hence $c_0\approx3\times10^{6}$ and $C_*\approx4\times10^{-8}$; \Cref{cor:supnorm} below then requires $\|v\|_{L^\infty}^2<\nu\|\omega\|_{L^\infty}/(9.6\times10^{7})$. The constant is dominated by the bound of $\nabla v$ in BMO and is far from optimal.
\end{remark}

\begin{remark}[Cancellation in the pressure source]\label{rem:source}
By~\eqref{eq:split}, the pressure source is $$\operatorname{tr}((\nabla v)^2)=|S|^2-\frac12|\omega|^2=|\nabla v|^2-|\omega|^2.$$ \Cref{cor:pressure} bounds it by $|\nabla v|^2$ in absolute value and so ignores any cancellation between strain and rotation. Any such cancellation lowers $c_0$ and raises $C_*$ in proportion. Since every term in the argument is critical, the gain affects only the constant.
\end{remark}

\begin{remark}[Large local Reynolds number]\label{rem:quartic}\sloppy
Hypothesis~\eqref{eq:hypmax} asks for a small local Reynolds number where the velocity dominates. In the opposite regime a damping can still come from the vorticity, through the quantity $|v|^4+\nu^2|\omega|^2$, whose second term carries the full dissipation $|\nabla\omega|^2$ in its equation. Take $f=0$; with a force and $\delta<0.45$ below, the additional terms are of lower order in the maximum, as in the proof of \Cref{thm:max}. From~\eqref{eq:NS} and the vorticity equation,
\begin{gather*}
(\partial_t+v\cdot\nabla-\nu\Delta)|v|^4+2\nu\big|\nabla|v|^2\big|^2+4\nu|v|^2|\nabla v|^2=-4|v|^2\,v\cdot\nabla p,\\
(\partial_t+v\cdot\nabla-\nu\Delta)\big(\nu^2|\omega|^2\big)+2\nu^3|\nabla\omega|^2=2\nu^2\,\omega\cdot(\omega\cdot\nabla)v .
\end{gather*}
Fix $t$ and let $x$ be a point where $|v|^4+\nu^2|\omega|^2$ attains its maximum $\Mq$. Its gradient vanishes there, so
\[
2\,|v|^3\,\big|\nabla|v|\big|=\nu^2\,|\omega|\,\big|\nabla|\omega|\big|\qquad\text{at }(x,t).
\]
Assume $v(x,t)\ne0$ and $\nabla|v|(x,t)\ne0$; then $\omega(x,t)\ne0$, because $\nabla|\omega|^2$ vanishes where $\omega$ does. The condition
\begin{equation}\label{eq:quartic-cond}
|\omega(x,t)|\le\big|\nabla|v|(x,t)\big|
\end{equation}
is then equivalent to $\nu^2|\nabla|\omega||\ge2|v|^3$ at $(x,t)$, and when it holds (Kato's inequality) that $2\nu^3|\nabla\omega|^2\ge8|v|^6/\nu$ there. Assume moreover that, for a constant $\delta\le0.45$,
\begin{equation}\label{eq:quartic-ball}
|\nabla v|\le\frac{\delta}{\nu}\,|v|^2\qquad\text{on the ball }B(x,r),\qquad r=\frac{\sqrt{11}\,\nu}{\delta\,\big(\Mq\big)^{1/4}} .
\end{equation}

The stretching is absorbed by the dissipation of $|v|^4$: by~\eqref{eq:quartic-ball} at $x$ and~\eqref{eq:quartic-cond}, $$2\nu^2\,\omega\cdot(\omega\cdot\nabla)v\le2\nu^2|\omega|^2|\nabla v|\le2\delta\nu|v|^2\big|\nabla|v|\big|^2,$$ while $2\nu\big|\nabla|v|^2\big|^2=8\nu|v|^2\big|\nabla|v|\big|^2$.
 For the pressure, split $\nabla p(x,t)=\nabla\Phi*\operatorname{tr}((\nabla v)^2)(x)$ at radius $r$ with a sharp cutoff. The component of $\nabla\Phi(z)$ along a unit vector $e$ is $-e\cdot z/(4\pi|z|^3)$, whose absolute value has integral $r/2$ over $B(0,r)$; since $|\operatorname{tr}((\nabla v)^2)|\le|\nabla v|^2\le\delta^2|v|^4/\nu^2$ on $B(x,r)$, the inner part is at most $$\frac{r\delta^2}{2\nu^2}\Mq.$$ In the outer part, $\operatorname{tr}((\nabla v)^2)=\nabla\cdot\nabla\cdot(v\otimes v),$ and two integrations by parts leave the integral over $|y-x|>r$ of the third derivatives of $\Phi$ against $v\otimes v$, and two integrals over the sphere $|y-x|=r$: of the second derivatives of $\Phi$ against $v\otimes v$, and of $\nabla\Phi$ against $(v\cdot\nabla)v$. Bounding each kernel by its largest eigenvalue and averaging over directions, and using $|v|^2\le(\Mq)^{1/2}$ everywhere, the first two are at most $$\big(3+\frac{3}{\sqrt5}\ln\frac{3+\sqrt5}{2}\big)(\Mq)^{1/2}/r<4.30\,(\Mq)^{1/2}/r$$ and $$\big(1+\frac{\ln(2+\sqrt3)}{4\sqrt3}\big)(\Mq)^{1/2}/r<1.20\,(\Mq)^{1/2}/r,$$ and the third is at most $\frac12\sup_{|y-x|=r}|v||\nabla v|\le\frac{\delta}{2\nu}(\Mq)^{3/4}$. With the value of $r$ in~\eqref{eq:quartic-ball},
\[
|\nabla p(x,t)|\le\Big(\sqrt{11}+\frac12\Big)\frac{\delta}{\nu}\,\big(\Mq\big)^{3/4} .
\]
Finally,~\eqref{eq:quartic-cond} and~\eqref{eq:quartic-ball} at $x$ give $\nu|\omega|\le\nu|\nabla v|\le\delta|v|^2$, so that $$\Mq\le(1+\delta^2)|v(x,t)|^4,$$ and $(\sqrt{11}+\frac12)\,\delta\,(1+\delta^2)^{3/4}\le2$ for $\delta\le0.45$. At $x$ the gradient of $|v|^4+\nu^2|\omega|^2$ vanishes and its Laplacian is nonpositive; adding the two equations and using the bounds above,
\begin{multline*}
\partial_t\big(|v|^4+\nu^2|\omega|^2\big)(x,t)\\
\le-\frac4\nu\,|v(x,t)|^3\Big(2|v(x,t)|^3-\Big(\sqrt{11}+\frac12\Big)\delta\,\big(\Mq\big)^{3/4}\Big)\le0 .
\end{multline*}
By the argument for maxima used in the proof of \Cref{thm:max}, the maximum does not increase at the times at which~\eqref{eq:quartic-cond} and~\eqref{eq:quartic-ball} hold at every maximum point. Nothing is gained where $\nabla|v|(x,t)=0$: there~\eqref{eq:quartic-cond} forces $\omega(x,t)=0$, and no damping is left.

Condition~\eqref{eq:quartic-cond} is restrictive. Together with~\eqref{eq:quartic-ball} it gives $\nu|\omega|\le\delta|v|^2$ at $(x,t)$, hence
\[
\frac{|\omega|}{\big|\nabla|\omega|\big|}\le\frac{\nu^2|\omega|}{2|v|^3}\le\frac\delta2\,\frac{\nu}{|v|}
\]
the vorticity must vary at $x$ on less than $\delta/2$ of the viscous length $\nu/|v|$. Solutions are analytic on that length: for data bounded by $M$, the solution extends analytically, up to a time of order $\nu/M^2$, to a complex strip of width of order $\sqrt{\nu t}$ on which it stays of order $M$~\cite{GK98}. By Cauchy's estimates, a time of order $\nu/M^2$ after an instant at which $|v|\le M$ everywhere, $|\nabla\omega|\le KM^3/\nu^2$ with $K$ an absolute constant. The damped case thus needs vorticity gradients within the factor $K/2$ of this bound, whereas the maximum may grow as soon as $\nu^2|\nabla|\omega||<2|v|^3$ at its maximum points.

The forced solution of~\cite{OpenAI26} is of the second kind. For $t$ close enough to $1$ it satisfies~\eqref{eq:quartic-ball} around the maximum points, where $|v|^2/|\nabla v|$ grows like $(1-t)^{-h}$ (\Cref{obs:II}) and where the radius in~\eqref{eq:quartic-ball}, of order $(1-t)^{1/2+h}$ with $\nu=1$, is small compared with the core radius $(1-t)^{1/2}$. In the core $|v|$ is of order $(1-t)^{-1/2-h}$ and $|\omega|$ of order $(1-t)^{-1-h}$. If $|\omega|$ varies on the core radius, then at a maximum point $|\nabla|v||/|\omega|=|\nabla|\omega||/(2|v|^3)$ is of order $(1-t)^{2h}$, and~\eqref{eq:quartic-cond} fails. Since $h<1/100$, this ordering appears only very close to $t=1$.
\end{remark}
\begin{corollary}\label{cor:bmo}
Let $c_0$ and $C_5$ be the constants of \Cref{cor:pressure} and \Cref{prop:bmo-interp}. If, for every $t\in(0,T^*)$,
\[
\|\nabla\times\NP(v(t))\|_{BMO}<C_5^{-1}\,\frac{\nu}{32c_0+1},
\]
then $T^*=T$.
\end{corollary}
The proof follows from the two lemmas below.
\begin{lemma}\label{prop:bmo-interp}
	There is an absolute constant $C_5$ such that, for every divergence-free $v\in H^3(\mathbb{R}^3)$,
	\[
	\|v\|_{L^\infty}^2\le C_5\,\|\nabla\times\NP(v)\|_{BMO}\,\|\omega\|_{L^\infty}.
	\]
\end{lemma}

\begin{proof}
	Since $-\Delta\NP$ is the identity and commutes with the curl, $-\Delta\nabla\times\NP(v)=\omega$, and since $\nabla\cdot v=0$, $\nabla\times\nabla\times\NP(v)=v$. For $\tau>0$,
	\[
	\nabla\times\NP(v)-e^{\tau\Delta}\,\nabla\times\NP(v)=-\int_0^\tau e^{s\Delta}\Delta\nabla\times\NP(v)\,ds=\int_0^\tau e^{s\Delta}\omega\,ds ,
	\]
	and taking the curl, which commutes with the heat semigroup,
	\[
	v(x)=\int_{\mathbb{R}^3}\nabla_x\Big((4\pi\tau)^{-\frac32}e^{-\frac{|x-y|^2}{4\tau}}\Big)\times\nabla\times\NP(v)(y)\,dy+\int_0^\tau\!\!\int_{\mathbb{R}^3}\nabla_x\Big((4\pi s)^{-\frac32}e^{-\frac{|x-y|^2}{4s}}\Big)\times\omega(y)\,dy\,ds .
	\]
	The gradient of the heat kernel has zero integral, so in the first integral $\nabla\times\NP(v)(y)$ may be replaced by its difference with its mean over $B(x,\sqrt\tau)$. On the annuli $2^k\sqrt\tau<|x-y|\le2^{k+1}\sqrt\tau$ the means differ from that one by at most $8(k+1)\|\nabla\times\NP(v)\|_{BMO}$, which the Gaussian decay absorbs, so the first integral is at most $C_4\tau^{-1/2}\|\nabla\times\NP(v)\|_{BMO}$ with $C_4$ absolute. Since $$\int\big|\nabla\big((4\pi s)^{-3/2}e^{-|z|^2/(4s)}\big)\big|\,dz=2/\sqrt{\pi s},$$ the second integral is at most $4\sqrt{\tau/\pi}\,\|\omega\|_{L^\infty}$. The choice $\tau=\frac{\sqrt\pi\,C_4}{4}\|\nabla\times\NP(v)\|_{BMO}/\|\omega\|_{L^\infty}$ makes the two bounds equal and gives the claim with $C_5=16C_4/\sqrt\pi$.
\end{proof}
\begin{lemma}\label{cor:supnorm}
	Let $c_0$ be the constant of \Cref{cor:pressure}. If, for every $t\in(0,T^*)$ at which $v(t)$ is not identically zero,
	\[
	\|v(t)\|_{L^\infty}^2<\frac{\nu}{32c_0+1}\,\|\omega(t)\|_{L^\infty},
	\]
	then $\sup_{t<T^*}\|v(t)\|_{L^\infty}<\infty$ and $T^*=T$.
\end{lemma}
\begin{proof}
	With $\lambda=\nu/(32c_0+1)$ as in \Cref{thm:max}, the maximum over $x$ of $\frac12|v|^2+\lambda(|\omega|^2+\varepsilon^2)^{1/2}$ at time $t$ is at least $\lambda\|\omega(t)\|_{L^\infty}$, which by hypothesis exceeds $|v(x,t)|^2$ at every point $x$. At every maximum point the first case of the proof of \Cref{thm:max} therefore occurs, and that case uses no hypothesis; the rest of that proof applies unchanged.
\end{proof}
\subsection{A sup-norm form and the Lamb vector}

The quantity of this subsection is the $L^1$ norm of the Lamb vector $v\times\omega$. Three facts single it out.

First, the Lamb vector is the whole nonlinearity up to a gradient. Since we have $$(v\cdot\nabla)v=\nabla\frac{|v|^2}{2}-v\times\omega,$$ the gradient part is absorbed by the pressure, and~\eqref{eq:NS} feels the nonlinearity only through the Leray projection of $v\times\omega$. The norm $\|v\times\omega\|_{L^1}$ therefore measures the effective nonlinearity itself, rather than the size of $v$ or of $\omega$. It is at most $\|v\|_{L^2}\|\omega\|_{L^2}=\|v\|_{L^2}\|\nabla v\|_{L^2}$, the quantity of the classical scale-invariant small-data condition, which goes back to~\cite{Leray34}, and it vanishes wherever $v$ and $\omega$ are parallel. Inequality~\eqref{eq:lamb-morrey} below makes this quantitative. A small Lamb vector forces a small kinetic energy in every ball, relative to its radius, and a divergence-free field in $H^3$ whose Lamb vector vanishes identically is zero: this is the Liouville property of Beltrami fields of finite energy~\cite{Nad14,CC15}.

Second, the norm is critical. $\|v\times\omega\|_{L^1}$ has the dimension of $\nu^2$ and is invariant under $v\mapsto a\,v(ax,a^2t)$. By \Cref{prop:lamb-bmo} it controls both the critical Morrey quantity $\sup_{x_0,r}r^{-1}\int_{B(x_0,r)}|v|^2$ and the norm in $BMO$ of the vector potential $\nabla\times\NP(v)$. Since $\nabla\cdot v=0$, $v=\nabla\times\big(\nabla\times\NP(v)\big)$, so every component of $v$ is a sum of first derivatives of components of a field of bounded mean oscillation. The second bound therefore places $v$ in $BMO^{-1}$, the space of first derivatives of functions of bounded mean oscillation. This space contains $L^3$, weak $L^3$ and the critical Besov spaces $\dot B^{-1+3/p}_{p,\infty}$, $p<\infty$, and it is the largest critical space in which small-data well-posedness is known~\cite{KT01}.

Third, this gives two regularity statements. \Cref{cor:lamb} goes through the maximum principle of \Cref{thm:max}: it allows a force, and asks for smallness at every time. Without a force, smallness at a single time suffices (\Cref{cor:lamb-single}), through the small-data theory in $BMO^{-1}$. Up to the values of the constants this is the stronger statement in the unforced case, and it contains the classical condition recalled above, with the product of the norms replaced by the norm of the product.

\begin{corollary}\label{cor:lamb}
	Let $c_0$ and $C_5$ be the constants of \Cref{cor:pressure} and \Cref{prop:bmo-interp}. If
	\[
	\sup_{0\le t<T^*}\|v(t)\times\omega(t)\|_{L^1}<\Big(\frac{\nu}{100\,C_5\,(32c_0+1)}\Big)^2,
	\]
	then $T^*=T$.
\end{corollary}
\begin{proposition}\label{prop:lamb-bmo}
For every divergence-free $v\in H^3(\mathbb{R}^3)$, with $\omega=\nabla\times v$,
\begin{equation}\label{eq:lamb-morrey}
\int_{B(x_0,r)}|v|^2\le16\,r\,\|v\times\omega\|_{L^1}\qquad\text{for all }x_0\in\mathbb{R}^3\text{ and }r>0,
\end{equation}
and
\[
\|\nabla\times\NP(v)\|_{BMO}\le\frac{304}{\sqrt{3\pi}}\,\|v\times\omega\|_{L^1}^{1/2}<100\,\|v\times\omega\|_{L^1}^{1/2},
\]
where $\|g\|_{BMO}$ is the supremum over balls $B$ of the average over $B$ of $\big|g-\frac1{|B|}\int_Bg\big|$.
\end{proposition}

\begin{proof}
It suffices to treat smooth compactly supported fields, which are dense among divergence-free fields in $H^3$, since all the quantities involved pass to the limit along such an approximation. Since $(v\cdot\nabla)v=\nabla{|v|^2}/{2}-v\times\omega$, we have \[\nabla\cdot(v\times\omega)=\Delta\frac{|v|^2}{2}-\nabla\cdot((v\cdot\nabla)v)=\Delta\frac{|v|^2}{2}-\nabla\cdot\nabla \cdot(v\otimes v).\] Convolving this identity with $|z|/(8\pi)$\footnote{ Which is $-\NP\circ\NP$ because $\Delta|z|=2/|z|.$} yields, \begin{enumerate}\item for the left-hand side, after one integration by parts, \[\frac1{8\pi}\int\frac{(x-y)\cdot(v\times\omega)(y)}{|x-y|}\,dy.\] \item for the right-hand one, \begin{enumerate} \item the first term being exactly $$\frac{1}{8\pi}\int\frac{|v(y)|^2}{|x-y|}\,dy;$$ 
		\item the second, integrated by parts twice against the Hessian of $|z|$\footnote{\(\nabla\nabla |z|=\frac{1}{|z|}\big(I-\frac{z\otimes z}{|z|^2}\big)\)}, giving, with the identity $|z|^2|v|^2=(z\cdot v)^2+|z\times v|^2$, \[-\frac1{8\pi}\int\frac{|(x-y)\times v(y)|^2}{|x-y|^3}\,dy=\frac1{8\pi}\int\frac{((x-y)\cdotp v(y))^2}{|x-y|^3}\,dy-\frac1{8\pi}\int\frac{|(x-y)|^2 |v(y)|^2}{|x-y|^3}\,dy\] \end{enumerate} after a sum up,  
		\[
		\int_{\mathbb{R}^3}\frac{\big((x-y)\cdot v(y)\big)^2}{|x-y|^3}\,dy=\int_{\mathbb{R}^3}\frac{(x-y)\cdot(v\times\omega)(y)}{|x-y|}\,dy\le\|v\times\omega\|_{L^1}.
		\]\end{enumerate}
Fix $B(x_0,r)$, average this over $x\in B(x_0,2r)$, and keep only $y\in B(x_0,r)$, for which $B(y,r)\subset B(x_0,2r)$. As we have \(\big((x-y)\cdot v(y)\big)^2=\rho^2|v(y)|^2\cos^2\theta\), polar coordinates around $y$ give \[\int_{B(y,r)}\frac{((x-y)\cdot v(y))^2}{|x-y|^3}\,dx=\frac{2\pi}{3}r^2|v(y)|^2,\] and dividing by $|B(x_0,2r)|={32\pi r^3}/{3}$ gives~\eqref{eq:lamb-morrey}.

For the second inequality, write $$\nabla\times\NP(v)(x)=\int\frac{v(y)\times(x-y)}{4\pi|x-y|^3}\,dy$$ and split the integral into the parts with $y\in B(x_0,2r)$ and $y\notin B(x_0,2r)$. For $y\in B(x_0,2r)$ the average of $|x-y|^{-2}$ over $x\in B(x_0,r)$ is at most $3/r^2$, its value for the ball of radius $r$ centered at $y$; by Fubini's theorem, the Cauchy-Schwarz inequality and~\eqref{eq:lamb-morrey} at radius $2r$, the average of the near part over $B(x_0,r)$ is at most \[\frac{3}{4\pi r^2}\int_{B(x_0,2r)}|v|\le\frac{24}{\sqrt{3\pi}}\|v\times\omega\|_{L^1}^{1/2}.\] The kernel \({a\times z}/{(4\pi|z|^3)}\) has $z$-gradient of norm at most $2|a|/(4\pi|z|^3)$, and for $x\in B(x_0,r)$ and $|y-x_0|>2r$ the segment from $x_0-y$ to $x-y$ stays at distance at least $\frac12|y-x_0|$ from the origin; Taylor's formula with integral remainder therefore bounds the difference between the far part at $x$ and its value at $x_0$ by \[\frac{4r}{\pi}\int_{|y-x_0|>2r}\frac{|v(y)|}{|y-x_0|^3}\,dy.\] 
On the shell $2^kr<|y-x_0|\le2^{k+1}r$, the Cauchy-Schwarz inequality and~\eqref{eq:lamb-morrey} at radius $2^{k+1}r$ bound the contribution by $2^{-k}128\|v\times\omega\|_{L^1}^{1/2}/\sqrt{3\pi}$, and summing over $k\ge1$ gives at most $$\frac{128}{\sqrt{3\pi}}\|v\times\omega\|_{L^1}^{1/2}.$$ 
The far part at $x_0$ is finite because $v\in L^2$, so the average over $B(x_0,r)$ of the distance between $\nabla\times\NP(v)$ and this constant is at most $\frac{152}{\sqrt{3\pi}}\|v\times\omega\|_{L^1}^{1/2}$; replacing the constant by the mean over the ball costs at most a factor $2$.
\end{proof}

\begin{remark}[Lamb's formula and the pressure]\label{rem:lamb-formula}
The identity in the proof of \Cref{prop:lamb-bmo} is a localised form of a classical one. Convolving the same identity with $|z|^2$ instead of $|z|$, whose Laplacian is $6$ and whose Hessian is $2I$, the same two integrations by parts give $2\int(x-y)\cdot(v\times\omega)(y)\,dy=3\int|v|^2-2\int|v|^2$, that is,
\[
\int_{\mathbb{R}^3}(x-y)\cdot(v\times\omega)(y)\,dy=\frac12\int_{\mathbb{R}^3}|v|^2\,dy ,
\]
Lamb's formula for the kinetic energy as the first moment of the Lamb vector~\cite{Lamb32,Saf92}. The kernel $|z|$ localises it: its Hessian $\frac{1}{|z|}\big(I-\frac{z\otimes z}{|z|^2}\big)$ removes the transverse part of the kinetic energy and leaves the radial kinetic energy seen from $x$, weighted by $1/|x-y|$. Equivalently, let $p$ be the pressure determined by $-\Delta p=\nabla\cdot\nabla\cdot(v\otimes v)$, so that $\Delta\big(p+\frac12|v|^2\big)=\nabla\cdot(v\times\omega)$. Then, for every $x$,
\[
\NP(p)(x)=-\frac1{8\pi}\int_{\mathbb{R}^3}\frac{|(x-y)\times v(y)|^2}{|x-y|^3}\,dy\le0,\quad
\NP\big(p+\tfrac12|v|^2\big)(x)=\frac1{8\pi}\int_{\mathbb{R}^3}\frac{\big((x-y)\cdot v(y)\big)^2}{|x-y|^3}\,dy\ge0 .
\]
The first formula is the analogue, one order smoother, of the classical singular-integral representation of the pressure, in which the same radial and transverse parts of the kinetic energy appear; see the local pressure formulae of~\cite{Con14}. We have not found in the literature the bound~\eqref{eq:lamb-morrey} of the critical Morrey quantity by the $L^1$ norm of the Lamb vector.
\end{remark}

\begin{proof}(of \Cref{cor:lamb})
At every $t$ at which $v(t)$ is not identically zero, \Cref{prop:lamb-bmo} and~\Cref{prop:bmo-interp} give
\[
\|v(t)\|_{L^\infty}^2\le100\,C_5\,\|v(t)\times\omega(t)\|_{L^1}^{1/2}\,\|\omega(t)\|_{L^\infty}<\frac{\nu}{32c_0+1}\,\|\omega(t)\|_{L^\infty},
\]
and \Cref{cor:supnorm} applies.
\end{proof}

\begin{corollary}\label{cor:lamb-single}
There is an absolute constant $\varepsilon_0>0$ with the following property. Let $f=0$, and assume that at a single time $t_0\in[0,T^*)$
\[
\|v(t_0)\times\omega(t_0)\|_{L^1}<\varepsilon_0\,\nu^2 .
\]
Then $T^*=T$ for every $T>t_0$: the solution is smooth on $[0,\infty)$.
\end{corollary}

\begin{proof}
By \Cref{prop:lamb-bmo}, $\|\nabla\times\NP(v(t_0))\|_{BMO}<100\,\varepsilon_0^{1/2}\,\nu$. Since $v(t_0)=\nabla\times\big(\nabla\times\NP(v(t_0))\big)$, every component of $v(t_0)$ is a sum of two first derivatives of components of $\nabla\times\NP(v(t_0))$, and a first derivative of a function of bounded mean oscillation belongs to $BMO^{-1}$, with norm at most an absolute constant $C_6$ times its norm in $BMO$~\cite{KT01}. Hence $\|v(t_0)\|_{BMO^{-1}}<100\,C_6\,\varepsilon_0^{1/2}\,\nu$. The change of variables $v(x,t)\mapsto\nu^{-1}v(x,t/\nu)$ turns~\eqref{eq:NS} with $f=0$ into the same equations with viscosity $1$ and divides the norm in $BMO^{-1}$ by $\nu$. By~\cite{KT01} there is therefore an absolute constant $\varepsilon_1>0$ such that, if $100\,C_6\,\varepsilon_0^{1/2}\le\varepsilon_1$, the datum $v(t_0)$ generates a global solution $u$ on $[t_0,\infty)$, smooth for $t>t_0$, with $\|u(t)\|_{L^\infty}\le C\,\nu^{1/2}(t-t_0)^{-1/2}$ for $t>t_0$ and an absolute constant $C$; we fix $\varepsilon_0$ accordingly. The smooth solution $v$ coincides with $u$ on $[t_0,T^*)$. Indeed, $v$ is bounded on $[t_0,t_0+\delta]$ for $\delta>0$ small, so its norm in the path space of~\cite{KT01} on that interval is at most an absolute constant times $\delta^{1/2}\nu^{-1/2}\sup_{t_0\le t\le t_0+\delta}\|v(t)\|_{L^\infty}$, which is small; both solutions then lie in the ball of that space in which the solution is unique, so they agree on $[t_0,t_0+\delta]$, and after $t_0+\delta$ they are smooth solutions with the same data. Hence $\|v(t)\|_{L^\infty}\le C\,\nu^{1/2}(t-t_0)^{-1/2}$ on $(t_0,T^*)$, which stays bounded as $t\to T^*$. A velocity bounded in $L^\infty$ cannot blow up, so $T^*=T$.
\end{proof}

With a force, the small-data theory in $BMO^{-1}$ also requires the force to be small after $t_0$ in a critical norm; we do not pursue this.

\begin{remark}\label{rem:lamb-scaling}
The smallness in \Cref{cor:lamb} is natural on dimensional and scaling grounds: $\|v\times\omega\|_{L^1}$ has the dimension of $\nu^2$ and is invariant under $v\mapsto a\,v(ax,a^2t)$, as are both propositions. This route therefore cannot remove the smallness; it can only move the threshold through the constants $c_0$, $100$ and $C_5$. Inequality~\eqref{eq:lamb-morrey} also shows that $\|v\times\omega\|_{L^1}$ controls the critical Morrey quantity $\sup_{x_0,r}r^{-1}\int_{B(x_0,r)}|v|^2$.
The same holds for \Cref{cor:lamb-single}, whose threshold $\varepsilon_0\nu^2$ comes from the small-data theory in $BMO^{-1}$.
\end{remark}

\begin{remark}[The Lamb vector along the forced blow-up]\label{rem:lamb-forced}
Along the forced solution of~\cite{OpenAI26} the norm $\|v\times\omega\|_{L^1}$ is unbounded, but it grows only through the anomalous factor. By \Cref{prop:bmo-interp} and \Cref{prop:lamb-bmo}, at every time
\[
\|v\|_{L^\infty}^2\le100\,C_5\,\|v\times\omega\|_{L^1}^{1/2}\,\|\omega\|_{L^\infty},
\]
and along that solution $\|v\|_{L^\infty}^2/\|\omega\|_{L^\infty}$ grows like $(1-t)^{-h}$ up to logarithmic factors (see the discussion after \Cref{obs:I}), so $\|v\times\omega\|_{L^1}$ grows at least like $(1-t)^{-2h}$. In the other direction, with the scalings recalled in \Cref{sec:construction} of the present paper, at height $z$ the radial integral of $|v||\omega|$ is of order $q^{-1/2-2h}$, where $q$ is the concentration scale; since $1+4h>1-2h$, its integral in $z$ is dominated by the core, of volume comparable to $(1-t)^{3/2-h}$, and $\int|v||\omega|\asymp(1-t)^{-3h}$, up to the logarithmic factors of the pulse annulus. Hence, at leading order,
\[
(1-t)^{-2h}\lesssim\|v\times\omega\|_{L^1}\lesssim(1-t)^{-3h}.
\]
Like the local Reynolds numbers of Criteria~\ref{crit:I} and~\ref{crit:II}, this critical quantity diverges only through the anomalous factor: for $h=1/100$ and $1-t=10^{-15}$, $(1-t)^{-3h}$ is below $3$. \Cref{cor:lamb-single} does not apply to this solution, since its force acts up to the singular time.
\end{remark}

\subsection{The signed pressure work}

\begin{proposition}\label{prop:signed}
Let $0\le\lambda<\sqrt3\,\nu$, $0\le\theta\le1-\lambda/(\sqrt3\,\nu)$, $A\ge0$ and $\eta\in(0,1)$. Assume that, for every $t\in(0,T^*)$,
\begin{equation}\label{eq:hypsigned}
-v(x,t)\cdot\nabla p(x,t)\le\theta\,\nu\,|\nabla v(x,t)|^2+A
\end{equation}
at every point $x$ where $m(x,t)\ge(1-\eta)\sup_ym(y,t)$. Then $\sup_{t<T^*}\|v(t)\|_{L^\infty}<\infty$ and $T^*=T$.
\end{proposition}

\begin{proof}
When $\lambda=0$ no regularisation is needed. When $\lambda>0$, a maximum point $x$ of $\frac12|v|^2+\lambda(|\omega|^2+\varepsilon^2)^{1/2}$ satisfies $m(x,t)\ge(1-\eta)\sup_ym(y,t)$ as soon as the maximum is at least $\lambda\varepsilon(1+1/\eta)$, as in the proof of \Cref{thm:max}. The bound~\eqref{eq:pointwise} at $x$ together with~\eqref{eq:hypsigned} gives
\begin{multline*}
D^+\Mx\\
\le A+\|f\|_{L^\infty}\Big(2\Mx\Big)^{1/2}+\lambda\|\nabla\times f\|_{L^\infty},
\end{multline*}
so the maximum grows at most quadratically in time and stays bounded on $[0,T]$. Letting $\varepsilon\to0$ bounds $\|v\|_{L^\infty}$ on $[0,T^*)$.
\end{proof}

\Cref{prop:signed} is elementary. For $\lambda=0$ it states that the maximum of $|v|$ can grow only where the pressure work $-v\cdot\nabla p$ exceeds the dissipation $\nu|\nabla v|^2$, the observation behind~\cite{TY16}. We include it because \Cref{sec:signed} shows that the forced blow-up of~\cite{OpenAI26} sits exactly at this balance.

\section{The criteria and the forced blow-up}\label{sec:forced}

We record four criteria and check them against the forced blow-up of~\cite{OpenAI26}. The statements of this section are leading-order consequences of the construction and of the scalings stated in~\cite{OpenAI26}; they are not complete proofs. Where the construction gives a scaling only up to logarithmic factors (the oscillatory pulses), our statements inherit those factors.

\subsection{The construction}\label{sec:construction}

By~\cite[Theorem~1.1]{OpenAI26}, for every $\nu>0$ there are a force $f\in C^\infty_c(\mathbb{R}^3\times(0,\infty))$, a compact set $K\subset\mathbb{R}^3$, and smooth $v,p$ on $\mathbb{R}^3\times[0,1)$ solving~\eqref{eq:NS} with $v(\cdot,0)=0$, with spatial support in $K$, with $\sup_{t<1}\|v(t)\|_{L^2}<\infty$ and with $\limsup_{t\to1}\|v(t)\|_{L^\infty}=\infty$. On $[0,1)$ this is a smooth solution with zero initial data and a smooth, compactly supported force, so it satisfies every hypothesis on the data made here. The viscosity enters only through $v(x,t)\mapsto\sqrt\nu\,v(x/\sqrt\nu,t)$, which leaves all local Reynolds numbers unchanged, so we take $\nu=1$.

The leading flow is an axisymmetric vortex concentrating at the origin. Its radial extent is comparable to $(1-t)^{1/2}$, its axial extent to $(1-t)^{1/2-h}$ with a fixed $0<h<1/100$, its azimuthal and axial velocities to $(1-t)^{-1/2-h}$, and its radial velocity is $O((1-t)^{-1/2})$. In cylindrical coordinates $(r,\theta,z)$ the azimuthal and axial components are, at leading order, a concentration scale raised to the power $-\frac12-h$ times smooth profiles of two similarity variables. At height $z$ the concentration scale is comparable to $(1-t)+|z|^{1/(1/2-h)}$; on the axis, at fixed axial similarity variable $\eta_0\in(-1,1)$, it equals $(1-t)/(1-\eta_0^2)$, at the moving height $z=\big((1-t)/(1-\eta_0^2)\big)^{1/2-h}\eta_0$. We use the swirl profile $F$ and the axial profile $U$ of~\cite[(4.3)-(4.4)]{OpenAI26}: $F$ is smooth and strictly positive, and $U(0,\cdot)$ is not identically zero, since a nonzero axis datum is what gives the lower bound on the axial velocity. Near the axis the corrections to the leading field are smaller by a factor $(1-t)^{2h}$, and the oscillatory pulses live in an annulus away from the axis; their amplitude $(1-t)^{-1/2-h/2}$ is below the core velocity, and their gradient, amplitude over wavelength $(1-t)^{1/2+h/2}$, is of the order $(1-t)^{-1-h}$ of the core gradient~\cite[Section~3.3]{OpenAI26}. Outside the annulus the flow is the purely azimuthal heat exterior, whose swirl is $c_\infty(r^2/2)^{-1/2-h}$ times a profile bounded between positive constants on the region considered below~\cite[(4.29)]{OpenAI26}.

From the Cartesian form of the leading field~\cite[(4.5)]{OpenAI26}, at the point of the axis with axial similarity variable $\eta_0$,
\begin{align}
v&=\Big(\frac{1-t}{1-\eta_0^2}\Big)^{-\frac12-h}U(0,\eta_0)\,e_z\,\big(1+O((1-t)^{2h})\big),\label{eq:axis-v}\\
\omega&=2\Big(\frac{1-t}{1-\eta_0^2}\Big)^{-1-h}F(0,\eta_0)\,e_z\,\big(1+O((1-t)^{2h})\big),\label{eq:axis-w}
\end{align}
and the Jacobian matrix of $v$ there (rows: components of $v$; columns: derivatives in $x$, $y$, $z$) is, at leading order,
\[
\nabla v=\Big(\frac{1-t}{1-\eta_0^2}\Big)^{-1-h}F(0,\eta_0)\begin{pmatrix}0&-1&0\\1&0&0\\0&0&0\end{pmatrix}+\Big(\frac{1-t}{1-\eta_0^2}\Big)^{-1}\frac{v_0(0,\eta_0)}{2}\begin{pmatrix}1&0&0\\0&1&0\\0&0&-2\end{pmatrix},
\]
where $v_0(0,\eta_0)$ is the value on the axis of the reduced radial-flux profile of~\cite[(4.4)]{OpenAI26} (not the initial datum, which vanishes). The antisymmetric part is the rotation of the core; the symmetric part, produced by the radial inflow, is smaller by the factor $(1-t)^{h}$, and its last diagonal entry follows from incompressibility.

\subsection{Criterion I: local Reynolds number on the high-vorticity set}

\begin{criterion}\label{crit:I}
There are a constant $C>0$ and a measurable $M\ge0$ with $\int_0^{T^*}M^{6/5}\,dt<\infty$ such that $|v|^2\le C\nu|\omega|$ wherever $|\omega|>M(t)$.
\end{criterion}

\Cref{thm:main} proves regularity under Criterion~\ref{crit:I} when $C<4/(3\sqrt3)$.

\begin{observation}\label{obs:I}
For every $C>0$ and every $M$ with $\int_0^1M^{6/5}\,dt<\infty$, the forced solution violates Criterion~\ref{crit:I}.
\end{observation}

Indeed, choose $\eta_0$ with $U(0,\eta_0)\ne0$. By~\eqref{eq:axis-v}-\eqref{eq:axis-w}, at the corresponding point of the axis
\begin{equation}\label{eq:ratio-axis}
\frac{|v|^2}{|\omega|}=(1-t)^{-h}\,\frac{(1-\eta_0^2)^{h}\,U(0,\eta_0)^2}{2F(0,\eta_0)}\big(1+O((1-t)^{2h})\big)\longrightarrow+\infty ,
\end{equation}
while $|\omega|\ge c\,(1-t)^{-1-h}$ there. Since $\int_{t_0}^1(1-t)^{-\frac65(1+h)}\,dt=\infty$ and $\int_{t_0}^1M^{6/5}\,dt<\infty$, the set of $t\in(t_0,1)$ with $M(t)<c\,(1-t)^{-1-h}$ has positive measure for every $t_0<1$. For $t_0$ close to $1$, at those times the axis point belongs to $\{|\omega|>M(t)\}$ and, by continuity, so does a neighbourhood of it on which $|v|^2>C|\omega|$. The sup-norm version fails as well: along the growth path of~\cite[Theorem~3.1(iv)]{OpenAI26}, $\|v(t)\|_{L^\infty}\ge c\,(1-t)^{-1/2-h}$, while the vorticity of the core and of the pulses is $O((1-t)^{-1-h})$ up to logarithmic factors, so $\|v\|_{L^\infty}^2/\|\omega\|_{L^\infty}\gtrsim(1-t)^{-h}$ up to logarithmic factors.

\begin{remark}
The violation is only by the anomalous factor $(1-t)^{-h}$, which for $h=1/100$ equals about $1.41$ at $1-t=10^{-15}$~\cite{Dur26}. Any variant of Criterion~\ref{crit:I} whose constant is allowed to grow in time as fast as $(1-t)^{-h}$ must be checked against this example.
\end{remark}

\subsection{Criterion II: local Reynolds number at the maximum points}

\begin{criterion}\label{crit:II}
For some $\lambda>0$, $\eta\in(0,1)$ and $C>0$, hypothesis~\eqref{eq:hypmax} holds at every point where $m\ge(1-\eta)\sup m$ and $|v|^2\ge(1-\eta)\sup m$.
\end{criterion}

\Cref{thm:max} proves regularity under Criterion~\ref{crit:II} when $\lambda=\nu/(32c_0+1)$ and $C<C_*$. Since $|\omega|\le\sqrt2\,|\nabla v|$, a pointwise bound $|v|^2\le C\nu|\omega|$ implies $|v|^2\le\sqrt2\,C\nu|\nabla v|$; the inequality of Criterion~\ref{crit:II} is the weaker of the two.

\begin{observation}\label{obs:II}
For every $\lambda>0$, $\eta\in(0,1)$ and $C>0$, the forced solution violates Criterion~\ref{crit:II}. Near the singular time the ratio $|v|^2/|\nabla v|$ is of order exactly $(1-t)^{-h}$ in the inner core, up to logarithmic factors in the pulse annulus.
\end{observation}

In the core $|v|^2$ is of order $(1-t)^{-1-2h}$ while $\lambda|\omega|$ is of order $(1-t)^{-1-h}$, so for $t$ close to $1$ the near-maximum points of $m$ lie where the velocity part dominates, and the hypothesis is required there. On the axis, the rotation part of the Jacobian dominates, $|\nabla v|=\sqrt2\,F(0,\eta_0)\big(\frac{1-t}{1-\eta_0^2}\big)^{-1-h}\big(1+O((1-t)^{2h})\big)$, and
\[
\frac{|v|^2}{|\nabla v|}=(1-t)^{-h}\,\frac{(1-\eta_0^2)^{h}\,U(0,\eta_0)^2}{\sqrt2\,F(0,\eta_0)}\big(1+O((1-t)^{2h})\big).
\]
Off the axis, on the plane $z=0$ at radius $r$ inside the inner core, the leading terms of $\nabla v$ are $u_\theta/r$, $\partial_ru_\theta$ and $\partial_ru_z$, all of order $(1-t)^{-1-h}$, and
\[
\frac{|v|^2}{|\nabla v|}\simeq(1-t)^{-h}\,\frac{E^2+U^2}{\Big(E^2\,\dfrac{1-t}{r^2}+\dfrac{r^2}{1-t}\big((\partial_1E)^2+(\partial_1U)^2\big)\Big)^{1/2}},
\]
where $E=\sqrt{2X}\,F$ is the swirl profile of~\cite{OpenAI26}, and $E$, $U$ and their derivatives $\partial_1E$, $\partial_1U$ in the first (radial similarity) variable are evaluated at $\big(r^2/(2(1-t)),0\big)$. The first term of the denominator equals $F^2>0$; it comes from $u_\theta/r$. The denominator therefore never vanishes in the core, and the ratio is bounded above and below by constants times $(1-t)^{-h}$ wherever $v\ne0$; at the near-maximum points of $|v|$ the numerator is bounded below, so the ratio tends to infinity there. In the pulse annulus the gradient is of order $(1-t)^{-1-h}$ up to logarithmic factors, so the ratio is $(1-t)^{-h}$ up to such factors, which $(1-t)^{-h}$ dominates. In the heat exterior, where the swirl is close to $c_\infty(r^2/2)^{-1/2-h}$ and $U=0$, the ratio is approximately $c_\infty(r^2/2)^{-h}$: independent of time at fixed $r$, and of order $(1-t)^{-h}$ at the inner edge $r\sim(1-t)^{1/2}$.

The flow offers no cancellation of the kind considered in \Cref{rem:source}. Near the axis the core is rotation-dominated, the pressure source is close to $-\frac12|\omega|^2$, and the pressure has its minimum on the axis, in centripetal balance $\partial_rp=u_\theta^2/r$. Off the axis~\cite{OpenAI26} gives $\big(\partial_ru_\theta-u_\theta/r,\ \partial_ru_z\big)=-q^{-A-1/2}F\,(a,-b_s)$, while the vertical vorticity is $\partial_ru_\theta+u_\theta/r$: strain and vorticity are both $(1-t)^{-1-h}$ times profile ratios that do not change in time. An improvement of $c_0$ through strain-rotation cancellation therefore cannot bring this flow inside Criterion~\ref{crit:II}; what diverges is the local Reynolds number itself.

The paired norms behave in the same way. At height $z$ the radial integral of a power of the fields contributes a factor equal to the concentration scale, and for the integrals below the integral in $z$ is dominated by $|z|\lesssim(1-t)^{1/2-h}$; the dominant region is the core, of volume comparable to $(1-t)^{3/2-h}$. Hence
\[
\int|v|^6\asymp(1-t)^{-\frac32-7h},\qquad\int|\omega|^3\asymp(1-t)^{-\frac32-4h},\qquad\int|\omega|^2|v|^2\asymp(1-t)^{-\frac32-5h},
\]
and $\|v\|_{L^6}^2/\|\omega\|_{L^3}\asymp(1-t)^{-h}$, a ratio invariant under the Navier-Stokes scaling. The three exponents realise equality in H\"older's inequality $\int|\omega|^2|v|^2\le(\int|\omega|^3)^{2/3}(\int|v|^6)^{1/3}$, and $\int_0^1\!\int|\omega|^3\,dx\,dt=\infty$. For comparison, the energy stays bounded and $\int|\omega|^2\asymp(1-t)^{-1/2-3h}$ is integrable in time, in agreement with the energy inequality.

\subsection{Criterion III: the vector potential}

Let $\NP$ be the inverse of $-\Delta$ on $\mathbb{R}^3$, $\NP(g)(x)=\frac{1}{4\pi}\int\frac{g(y)}{|x-y|}\,dy$. For divergence-free $v$, $\nabla\times\NP(v)=\NP(\omega)$ and $\nabla\times\nabla\times\NP(v)=v$, so $\nabla\times\NP(v)$ is the divergence-free vector potential of $v$, and
\[
\nabla\times\NP(v)(x)=\frac{1}{4\pi}\int\frac{y-x}{|y-x|^3}\times v(y)\,dy .
\]

\begin{proposition}\label{prop:interp}
For every $v\in L^2\cap L^\infty(\mathbb{R}^3)$,
\[
\|\nabla\times\NP(v)\|_{L^\infty}\le\Big(1+\frac{1}{\sqrt{4\pi}}\Big)\|v\|_{L^2}^{2/3}\,\|v\|_{L^\infty}^{1/3}.
\]
The exponents are sharp for general fields.
\end{proposition}

\begin{proof}
Bounding the kernel by $|x-y|^{-2}$ and splitting at distance $R$ from $x$, the ball contributes at most $\frac{1}{4\pi}\|v\|_{L^\infty}\int_{|x-y|<R}|x-y|^{-2}dy=R\|v\|_{L^\infty}$, and by the Cauchy-Schwarz inequality the complement contributes at most $\frac{1}{4\pi}\|v\|_{L^2}\big(\int_{|x-y|>R}|x-y|^{-4}dy\big)^{1/2}=\|v\|_{L^2}/\sqrt{4\pi R}$. Choosing $R=(\|v\|_{L^2}/\|v\|_{L^\infty})^{2/3}$ gives the bound. A field equal to a constant vector of size $m$ on a ball, with the radius chosen so that the energy is fixed, shows that the exponents cannot be improved.
\end{proof}

\begin{criterion}\label{crit:III}
Near the singular time the bound of \Cref{prop:interp} is reversed: there is a constant $K$ such that $\|v\|_{L^2}^{2/3}\|v\|_{L^\infty}^{1/3}\le K\|\nabla\times\NP(v)\|_{L^\infty}$. Since \Cref{prop:interp} always holds, this amounts to the equivalence of the two sides.
\end{criterion}

\begin{observation}\label{obs:III}
For the forced solution,
\[
\|\nabla\times\NP(v)(t)\|_{L^\infty}\asymp h^{-1}(1-t)^{-h}.
\]
Hence $\|v\|_{L^2}^{2/3}\|v\|_{L^\infty}^{1/3}\big/\|\nabla\times\NP(v)\|_{L^\infty}\asymp h\,(1-t)^{-1/6+2h/3}\to\infty$, and Criterion~\ref{crit:III} is violated by a positive power of $1-t$, not only by the anomalous factor.
\end{observation}

At the origin, $y\times e_\theta=r\,e_z-z\,e_r$ and the radial part averages to zero in the angle, while the poloidal components of $v$ produce only azimuthal vectors, which also average to zero. Thus $\nabla\times\NP(v)(0,t)$ is vertical and
\[
e_z\cdot\nabla\times\NP(v)(0,t)=\frac12\int_{-\infty}^{\infty}\!\int_0^\infty\frac{r^2\,u_\theta(r,z,t)}{(r^2+z^2)^{3/2}}\,dr\,dz ,
\]
with a positive integrand, since the swirl profile is positive off the axis. On the region $|z|\le r$, $C(1-t)^{1/2}\le r\le r_0$ the concentration scale is at most a constant times $(1-t)+r^{2/(1-2h)}$, which is small compared with $r^2$, so the region lies in the heat exterior, where $u_\theta\ge c\,r^{-1-2h}$. Since $\int_{|z|\le r}(r^2+z^2)^{-3/2}dz\ge1/(\sqrt2\,r^2)$,
\[
e_z\cdot\nabla\times\NP(v)(0,t)\ge\frac{c}{2\sqrt2}\int_{C(1-t)^{1/2}}^{r_0}r^{-1-2h}\,dr=\frac{c}{2\sqrt2}\,\frac{\big(C(1-t)^{1/2}\big)^{-2h}-r_0^{-2h}}{2h}.
\]
For the upper bound, near the singular point $|v(y)|\lesssim\min\{(1-t)^{-1/2-h},\,r(y)^{-1-2h}\}$ up to bounded terms, with $r(y)$ the distance to the axis. Integrating the kernel $|x-y|^{-2}$ first along the axis direction reduces the bound to a planar integral of $\min\{(1-t)^{-1/2-h},r^{-1-2h}\}$ against $\pi/|x_h-y_h|$, which is at most a constant times $h^{-1}(1-t)^{-h}$ plus a bounded term; the pulses contribute amplitude times wavelength, which is $O(1)$ up to logarithmic factors. Since $\|v\|_{L^\infty}\asymp(1-t)^{-1/2-h}$, equivalently $\|\nabla\times\NP(v)\|_{L^\infty}\asymp\|v\|_{L^\infty}^{2h/(1+2h)}$, an exponent below $0.02$. The interpolation of \Cref{prop:interp} is lossy here because saturating it would require $|v|$ close to $\|v\|_{L^\infty}$ on a ball of radius $(\|v\|_{L^2}/\|v\|_{L^\infty})^{2/3}\asymp(1-t)^{1/3+2h/3}$, whereas the core is a slender column of radius $(1-t)^{1/2}$ whose energy tends to zero.

Two further remarks locate the mechanism. First, $\nabla\times\NP(v)$ solves
\[
\partial_t\,\nabla\times\NP(v)-\nu\Delta\nabla\times\NP(v)=\nabla\times\NP(v\times\omega)+\nabla\times\NP(f),
\]
with zero initial data. The force term obeys $\|\nabla\times\NP(f)\|_{L^\infty}\le\|f\|_{L^\infty}+\|f\|_{L^1}/(4\pi)$, which stays bounded, and the heat semigroup does not increase the sup norm; by Duhamel's formula the whole divergence of $\|\nabla\times\NP(v)\|_{L^\infty}$ is produced by the Lamb-vector source $\nabla\times\NP(v\times\omega)$. Second, the reverse inequality of Criterion~\ref{crit:III} is false for general fields, through rapid oscillation or concentration into slender columns; if it holds, it holds because of the dynamics.

A bound on $\|\nabla\times\NP(v)\|_{L^\infty}$ alone would be an endpoint criterion. Since $\|v\|_{\dot B^{-1}_{\infty,\infty}}\le C\|\nabla\times\NP(v)\|_{L^\infty}$, it places $v$ in $L^\infty(0,T;\dot B^{-1}_{\infty,\infty})$, the largest critical space, where regularity is proved under a smallness condition on the jumps in time, in particular under continuity in time~\cite{CS10}, and is claimed under boundedness alone in a recent preprint~\cite{Ri23}. Since $|\nabla\times\NP(v)|=|\NP(\omega)|\le\NP(|\omega|)\le\sqrt2\,\NP(|\nabla v|)$, such a criterion would also give regularity whenever $\NP(|\nabla v|)$ stays bounded. The estimates of this paper give no criterion at such an endpoint; that would require a limiting argument of the kind behind the $L^\infty(0,T;L^3)$ endpoint~\cite{ESS03}, which we leave for later work.

\begin{remark}[Closure by Gronwall under a power bound on $\|v\|_{L^\infty}$]\label{rem:gronwall-power}
Let $v$ be a smooth solution of~\eqref{eq:NS} with $\nu=1$ on $[0,T)$ and let $L$ be its vector potential, $\nabla\times L=v$, $\nabla\cdot L=0$, so that $L=\nabla\times\NP(v)$ and, applying $\nabla\times\NP$ to the equation, $\partial_tL-\Delta L+\nabla\times\NP(\omega\times v)=\nabla\times\NP(f)$. The nonlinear term obeys
\[
\|\nabla\times\NP(\omega\times v)\|_{L^\infty}\le C\,\|v\|_{L^\infty}^{1/2}\,\|\Delta v\|_{L^2}^{1/2}\,\|\omega\times v\|_{L^2}^{1/2}.
\]
Indeed, writing $\nabla\times\NP(\omega\times v)(x_0)=\int_{\mathbb{R}^3}\nabla\Phi(x_0-y)\times(\omega\times v)(y)\,dy$ with $|\nabla\Phi(z)|=(4\pi|z|^2)^{-1}$ and splitting the integral at radius $r$, the inner part is bounded by H\"older's inequality and the Sobolev bound $\|\omega\|_{L^6}\le C\|\nabla\omega\|_{L^2}=C\|\Delta v\|_{L^2}$,
\[
\frac{\|v\|_{L^\infty}}{4\pi}\int_{B_r(x_0)}\frac{|\omega(y)|}{|x_0-y|^2}\,dy\le\frac{\|v\|_{L^\infty}}{4\pi}\,\|\omega\|_{L^6}\Big(\int_{B_r(0)}|y|^{-12/5}\,dy\Big)^{5/6}\le C\,r^{1/2}\,\|v\|_{L^\infty}\,\|\Delta v\|_{L^2},
\]
while the outer part is bounded by the Cauchy-Schwarz inequality against $\int_{|y|>r}|y|^{-4}\,dy=4\pi/r$, which gives $(4\pi r)^{-1/2}\|\omega\times v\|_{L^2}$; the choice $r=\|\omega\times v\|_{L^2}\big/\big(\|v\|_{L^\infty}\|\Delta v\|_{L^2}\big)$ balances the two contributions.

Suppose now that, for some $p\ge2$ and some constant $C$, $\|v(t)\|_{L^\infty}\le C\,\|L(t)\|_{L^\infty}^{p}$ on $[0,T)$. Only the factor $\|v\|_{L^\infty}^{1/p}$ is traded, $\|v\|_{L^\infty}^{1/2}=\|v\|_{L^\infty}^{\frac12-\frac1p}\|v\|_{L^\infty}^{\frac1p}\le C\,\|v\|_{L^\infty}^{\frac12-\frac1p}\,\|L\|_{L^\infty}$, so that the bound stays linear in $L$. Since the heat semigroup is a contraction on $L^\infty$,
\[
\|L(t)\|_{L^\infty}\le\|L(0)\|_{L^\infty}+\int_0^t\|\nabla\times\NP(f)\|_{L^\infty}\,ds+C\int_0^t\|L(s)\|_{L^\infty}\,\|v\|_{L^\infty}^{\frac12-\frac1p}\,\|\Delta v\|_{L^2}^{1/2}\,\|\omega\times v\|_{L^2}^{1/2}\,ds,
\]
and Gronwall's lemma yields
\[
\|L(t)\|_{L^\infty}\le\Big(\|L(0)\|_{L^\infty}+\int_0^t\|\nabla\times\NP(f)\|_{L^\infty}\,ds\Big)\exp\Big(C\int_0^t\|v\|_{L^\infty}^{\frac12-\frac1p}\,\|\Delta v\|_{L^2}^{1/2}\,\|\omega\times v\|_{L^2}^{1/2}\,ds\Big).
\]
Hence $L$ stays bounded on $[0,T)$ as soon as $\int_0^T\|v\|_{L^\infty}^{\frac12-\frac1p}\|\Delta v\|_{L^2}^{1/2}\|\omega\times v\|_{L^2}^{1/2}\,dt<\infty$; then $\|v\|_{L^\infty}\le C\|L\|_{L^\infty}^{p}$ stays bounded as well, and a velocity bounded in $L^\infty$ cannot blow up, so the solution continues past $T$.

Under the scaling $v_\lambda(x,t)=\lambda v(\lambda x,\lambda^2t)$ the time integral in the exponent scales like $\lambda^{-1/p}$: it is supercritical for every finite $p$ and becomes critical only in the limit $p\to\infty$. The substance of the remark therefore lies in the hypothesis $\|v\|_{L^\infty}\le C\|L\|_{L^\infty}^{p}$, which bounds $\|v\|_{L^\infty}$ by a power of the scale-invariant quantity $\|L\|_{L^\infty}$.

The forced blow-up of~\cite{OpenAI26} is consistent with this. At the leading order of its self-similar profile, written in the anisotropic similarity variables of the construction (see also~\cite{Dur26}) with anisotropy exponent $0<h<1/100$ and $\tau=T-t$, one has $\|v\|_{L^\infty}\simeq\tau^{-\frac12-h}$, $\|\Delta v\|_{L^2}\simeq\tau^{-\frac34-\frac{3h}{2}}$, $\|\omega\times v\|_{L^2}\simeq\tau^{-\frac34-\frac{5h}{2}}$ and $\|L\|_{L^\infty}\simeq\tau^{-h}$ (logarithmic in $\tau$ as long as $h\ln\frac1\tau\lesssim1$). The integrand in the exponent then behaves like $\tau^{-(\frac12+h)(\frac12-\frac1p)-\frac34-2h}$, which is integrable if and only if $p<\frac{1+2h}{5h}$, whereas $\|v\|_{L^\infty}\le C\|L\|_{L^\infty}^{p}$ holds along that solution if and only if $p\ge\frac{1+2h}{2h}$. The two ranges are disjoint, their thresholds differing by the factor $5/2$ whatever $h$ is, as they must be for a genuine blow-up.
\end{remark}

No regularity theorem under Criterion~\ref{crit:III} is proved here; the criterion is recorded for the margin by which the forced solution violates it.

\subsection{Criterion IV: the signed pressure work at the velocity maximum}\label{sec:signed}

\begin{criterion}\label{crit:IV}
Hypothesis~\eqref{eq:hypsigned} holds at the near-maximum points of $m$, with $\theta\le1-\lambda/(\sqrt3\,\nu)$.
\end{criterion}

\Cref{prop:signed} proves regularity under Criterion~\ref{crit:IV}.

\begin{observation}\label{obs:IV}
In the core of the forced solution, $-v\cdot\nabla p$ and $\nu|\nabla v|^2$ are both of order $(1-t)^{-2-2h}$, while $|v||\nabla p|$ is of order $(1-t)^{-2-3h}$. The solution violates Criterion~\ref{crit:IV}, necessarily so since $\|v\|_{L^\infty}\to\infty$, but by a bounded ratio and not by a power of $1-t$.
\end{observation}

In the core $\nabla p$ is essentially radial, given by the centripetal balance $\partial_rp=u_\theta^2/r$ of order $(1-t)^{-3/2-2h}$, with an axial component $\partial_zp$ of order $(1-t)^{-1-2h}/(1-t)^{1/2-h}=(1-t)^{-3/2-h}$, while $v$ is essentially azimuthal and axial. The signed pressure work is therefore
\[
-v\cdot\nabla p=-u_r\,\frac{u_\theta^2}{r}-u_z\,\partial_zp+\dots,
\]
with both displayed terms of order $(1-t)^{-1/2}(1-t)^{-3/2-2h}=(1-t)^{-1/2-h}(1-t)^{-3/2-h}=(1-t)^{-2-2h}$, the order of $\nu|\nabla v|^2$. The unsigned product $|v||\nabla p|$, of order $(1-t)^{-1/2-h}(1-t)^{-3/2-2h}=(1-t)^{-2-3h}$, loses exactly the anomalous factor through the angle between the swirl and the pressure gradient. In the central part of the core, where the fluid spirals inward, $-u_r\,u_\theta^2/r=|u_r|\,u_\theta^2/r>0$: fluid moving inward down the centripetal pressure gradient gains kinetic energy at the same order at which viscosity removes it. By \Cref{prop:signed} the work must exceed $\theta\nu|\nabla v|^2+A$ at near-maximum points for times accumulating at $1$, which is consistent with the growth of the leading velocity scale: its square $(1-t)^{-1-2h}$ has time derivative of order $(1-t)^{-2-2h}$, that is, of the order of its own power $1+1/(1+2h)$. This growth law is superlinear and reaches infinity in finite time.

\subsection{Summary}

\begin{center}
\renewcommand{\arraystretch}{1.35}
\begin{tabular}{@{}>{\raggedright}p{0.16\linewidth}>{\raggedright}p{0.31\linewidth}>{\raggedright}p{0.27\linewidth}>{\raggedright\arraybackslash}p{0.18\linewidth}@{}}
\toprule
Criterion & Hypothesis & Along the forced blow-up & Margin of violation\\
\midrule
I (\Cref{thm:main}) & $|v|^2\le C\nu|\omega|$ where $|\omega|>M(t)$ & $|v|^2/|\omega|\asymp(1-t)^{-h}$ on the axis & anomalous factor $(1-t)^{-h}$\\
II (\Cref{thm:max}) & $|v|^2\le C\nu|\nabla v|$ at near-maximum points of $m$ & $|v|^2/|\nabla v|\asymp(1-t)^{-h}$ in the core & anomalous factor $(1-t)^{-h}$\\
III (Vector potential) & $\|v\|_{L^2}^{2/3}\|v\|_{L^\infty}^{1/3}\le K\|\nabla\times\NP(v)\|_{L^\infty}$ & ratio $\asymp h\,(1-t)^{-1/6+2h/3}$ & power of $1-t$\\
IV (\Cref{prop:signed}) & $-v\cdot\nabla p\le\theta\nu|\nabla v|^2+A$ at near-maximum points of $m$ & $-v\cdot\nabla p$ and $\nu|\nabla v|^2$ both $\asymp(1-t)^{-2-2h}$ & bounded ratio\\
\bottomrule
\end{tabular}
\end{center}

In each case the forced solution of~\cite{OpenAI26} lies outside the hypothesis, so none of the four criteria contradicts the construction. Criteria I and II bound local Reynolds numbers and are violated only through the anomalous factor; Criterion III is violated robustly; Criterion IV is violated by a bounded ratio and sits at the balance the construction exploits. Because the construction uses a smooth compactly supported force, any proof that unforced Leray solutions~\cite{Leray34} satisfy one of these hypotheses must use the absence of the force, or a restriction on its structure, at some point: an argument valid for arbitrary smooth forces would also apply to this solution.

\section*{Declaration on the use of generative AI}

{
This work was carried out with the assistance of Claude (Anthropic; model Claude Opus 5.5), used through the Claude app in November 2026. The use went beyond language editing; the contributions are as follows.

\paragraph{Author.} The regularity criteria and the multipliers on which they rest are the author's: the insertion of a multiple of $|\omega|v$ into the pairing with $\omega\times v$ and the identity~\eqref{eq:identity}; the maximum principle for $\frac12|v|^2+|\omega|$ with the pressure controlled through the weighted estimate of \Cref{lem:weighted}; the multiplier $\NP(|\omega|^2)v$ of \Cref{thm:newton}, with the bound on $\int\NP(|\omega|^2)^2|v|^2$ through the $L^{6/5}$ norm of $|\omega|^2$; the gradient $\nabla\NP(|\omega|^2)$ added to $\omega\times v$ in \Cref{thm:coulomb}, with the reduction of its cross term to integrals of $\NP(|\omega|^2)$ against $|\omega|^2$ and $|\nabla v|^2$; and the observation of \Cref{rem:quartic}.

\paragraph{Assistant.} The assistant found the exact cancellation of the integrals of $|\omega|^2|v|^2$ in \Cref{thm:newton} and the integration by parts that avoids an $L^4$ norm of the vorticity. It wrote out the proofs of \Cref{thm:newton,thm:coulomb}, and the proof of \Cref{lem:weighted} at the author's direction. It carried out the bookkeeping of constants in \Cref{thm:main,thm:max}, including the improvement of the low-vorticity condition to $\int M^{6/5}<\infty$ and the stretching bound of \Cref{prop:pointwise}, and completed the estimates of \Cref{rem:quartic}. It proposed the estimate~\eqref{eq:lamb-morrey} of \Cref{prop:lamb-bmo}, in a discussion with the author of the $L^1$ norm of the Lamb vector and of the equation satisfied by $\nabla\times\NP(v)$, together with \Cref{cor:lamb-single} and the observation that~\eqref{eq:lamb-morrey} contains the Liouville property of Beltrami fields. It carried out the leading-order comparisons with the forced construction~\cite{OpenAI26} in \Cref{sec:forced}; these are consequences of the scalings stated in~\cite{OpenAI26}, not complete proofs. It also contributed the remarks of \Cref{sec:scaling} on changes of variables and \Cref{rem:M}, searched the literature and typeset the manuscript.

All mathematical arguments, computations and references have been checked by the author, who takes full responsibility for the content of the article. In accordance with the policies of arXiv and of the publisher, the assistant is not listed as an author.
}

\end{document}